%% file: Manuscript_updated.tex
\documentclass[preprint]{elsarticle} 

\makeatletter
\def\ps@pprintTitle{%
\let\@oddhead\@empty
\let\@evenhead\@empty
\def\@oddfoot{\centerline{\thepage}}%
\let\@evenfoot\@oddfoot}
\makeatother

\usepackage{etoolbox}
\patchcmd{\MaketitleBox}{\footnotesize\itshape\elsaddress\par\vskip36pt}{\footnotesize\itshape\elsaddress\par\parbox[b][36pt]{\linewidth}{\vfill\hfill\textnormal{\today}\hfill\null\vfill}}{}{}%
\patchcmd{\pprintMaketitle}{\footnotesize\itshape\elsaddress\par\vskip36pt}{\footnotesize\itshape\elsaddress\par\parbox[b][36pt]{\linewidth}{\vfill\hfill\textnormal{\today}\hfill\null\vfill}}{}{}%

\usepackage[
colorlinks=true,%
breaklinks=true,%
linkcolor=blue,
urlcolor=blue,%
citecolor=blue,%
pdftitle={Title of paper}, 
pdfkeywords={Keywords},	
pdfauthor={Authors},		
bookmarksopen=false,
pdfpagemode=None]{hyperref}
\usepackage{algorithm}
\usepackage{algorithmicx}
\usepackage{algpseudocode}
\usepackage{tabularx}
\usepackage{array}
\usepackage{multirow}
\usepackage{threeparttable}
\usepackage{booktabs}
\usepackage{graphicx}
\usepackage{arydshln}
\usepackage{setspace}   
\input{ListPackages}

\usepackage[textsize=tiny]{todonotes}

\newcommand{\RNum}[1]{\uppercase\expandafter{\romannumeral #1\relax}}
\usepackage{amsthm}
\makeatletter
\def\els@aparagraph[#1]#2{\elsparagraph[#1]{#2\@addpunct{.}}}
\def\els@bparagraph#1{\elsparagraph*{#1\@addpunct{.}}}
\makeatother
\input{mymacros.tex}

\begin{document}	
\begin{frontmatter}
		
\title{Multifidelity conditional value-at-risk estimation \\ by dimensionally decomposed generalized polynomial chaos-Kriging}
				
\author[ucsd]{Dongjin Lee\corref{cor1}}
\ead{dongjin-lee@ucsd.edu}		
\author[ucsd]{Boris Kramer\corref{cor1}}
\ead{bmkramer@ucsd.edu}
		
\cortext[cor1]{Corresponding author}
		
\address[ucsd]{Department of Mechanical and Aerospace Engineering, University of California San Diego, CA, United States}

\begin{abstract}
We propose novel methods for Conditional Value-at-Risk (CVaR) estimation for nonlinear systems under high-dimensional dependent random inputs. We develop a novel DD-GPCE-Kriging surrogate that merges dimensionally decomposed generalized polynomial chaos expansion and Kriging to accurately approximate nonlinear and nonsmooth random outputs. We use DD-GPCE-Kriging (1) for Monte Carlo simulation (MCS) and (2) within multifidelity importance sampling (MFIS). The MCS-based method samples from DD-GPCE-Kriging, which is efficient and accurate for high-dimensional dependent random inputs,  yet introduces bias. Thus, we propose an MFIS-based method where DD-GPCE-Kriging determines the biasing density, from which we draw a few high-fidelity samples to provide an unbiased CVaR estimate. To accelerate the biasing density construction, we compute DD-GPCE-Kriging using a cheap-to-evaluate low-fidelity model. Numerical results for mathematical functions show that the MFIS-based method is more accurate than the MCS-based method when the output is nonsmooth. The scalability of the proposed methods and their applicability to complex engineering problems are demonstrated on a two-dimensional composite laminate with 28 (partly dependent) random inputs and a three-dimensional composite T-joint with 20 (partly dependent) random inputs. In the former, the proposed MFIS-based method achieves $104$x speedup compared to standard MCS using the high-fidelity model, while accurately estimating CVaR with $1.15$\% error.
\end{abstract}

		
\end{frontmatter}

\section{Introduction} \label{sec:intro}
Risk measures are indispensable to designing safe and reliable systems in the presence of uncertainties that arise from many sources, such as manufacturing processes, material properties, or operating environments. In risk-averse design optimization \cite{Chaudhuri2022,wu2021risk,rocchetta2021scenario}, risk measures are used to avoid or handle a \emph{risky} output via objective and/or constraint functions. Here, we focus on two representative risk measures (originally employed in finance~\cite{RTRockafellar_SUryasev_2002a} but nowadays used in engineering fields): the Value-at-Risk (VaR) and the Conditional Value-at-Risk (CVaR). While the VaR aggregates losses (related to costs or target violation in any performance function) with a quantile of the distribution for the losses, the CVaR  reflects the mean or average value of losses exceeding VaR. The CVaR has superior quantitative or qualitative benefits over VaR, which has been demonstrated in portfolio management~\cite{rockafellar2000optimization,RTRockafellar_SUryasev_2002a,RMansini_WOgryczak_MGSperanza_2007a} and engineering design~\cite{HYang_MGunzburger_2016a,JORoyset_LBonfiglio_GVernengo_SBrizzolara_2017a}.
For example, CVaR (as a coherent risk measure) measures the average tail risk and preserves subadditivity and convexity; it thus facilitates optimization~\cite{rockafellar2000optimization,kouri2016risk,rocchetta2021scenario,Chaudhuri2022}. Nevertheless, optimization may require smoothed approximations of CVaR~\cite{kouri2016risk,kouri2020epi}. These benefits make CVaR appealing compared to other possibilities of measuring risks, such as reliability or failure probability, which is widely used in engineering~\cite{tyrrell2015engineering,Chaudhuri2022,wu2021risk}. Recently, $\mathrm{CVaR}$ has been incorporated into design of civil~\cite{tyrrell2015engineering}, energy~\cite{stover2023reliability}, naval~\cite{JORoyset_LBonfiglio_GVernengo_SBrizzolara_2017a}, and  aerospace~\cite{MANCUSO2021opt,JORoyset_LBonfiglio_GVernengo_SBrizzolara_2017a,chaudhuri2020multifidelity} applications.

CVaR and related risk measures are mostly estimated via sampling, e.g., Monte Carlo Simulation (MCS), especially when the quantity of interest does not follow a parametric distribution. Such sampling methods generally require a large number of samples to capture the tail risk of the quantity of interest. When each sample requires an expensive high-fidelity model run, e.g., a finite element analysis (FEA), standard MCS becomes computationally intensive, if not prohibitive.  
Numerous surrogate methods (e.g., reduced-order modeling (ROM) \cite{HKTQ18CVaRROMS,ZZou_DPKouri_WAquino_2019a,HKT2020_Adaptive_ROM_CVAR_estimation}, polynomial chaos expansion (PCE)~\cite{bernal2020volatility}, Kriging or Gaussian process modeling~\cite{chen2012stochastic}, support vector machine~\cite{gotoh2017support}, and neural networks~\cite{soma2020statistical}) have been used to break that computational bottleneck. Recently,~\cite{JAKEMAN2022108280} presented a surrogate model strategy to conservatively estimate risk measures from given data.
Many of these methods, especially the expansion or decomposition methods, are restricted to the assumption that the input random variables are independent.  
In reality, there is often correlation or dependence among some, or all, of the random inputs. 
Whether emanating from loads, material properties, or manufacturing variables, neglecting these correlations or dependencies may produce inefficient or risky designs~\cite{noh09,lee2020practical}.

Some PCE methods, e.g., the generalized PCE (GPCE)~\cite{rahman2018polynomial}, dimensionally decomposed GPCE (DD-GPCE)~\cite{lee2023high} or other PCE variants~\cite{Navarro2014POLYNOMIALCE,JAKEMAN2019643}, can handle dependent input random variables directly without a potentially detrimental measure transformation between dependent and independent random variables.  
A recent GPCE-based surrogate~\cite{lee2020practical} generates a multivariate orthonormal polynomial basis consistent with any non-product-type probability measure of inputs numerically instead of an analytical expression by a Rodrigues-type formula in~\cite{rahman2018polynomial}. DD-GPCE has been demonstrated to mitigate the curse of dimensionality to some extent by reshuffling and pruning the GPCE basis in a dimensionwise manner~\cite{lee2023high}.  In~\cite{lee2022bi}, we used DD-GPCE to estimate CVaR for a nonlinear output under high-dimensional and dependent random inputs and showed that DD-GPCE can be efficiently computed via a bifidelity method. However, DD-GPCE or other PCE-related methods generally have limitations when handling highly nonlinear or nonsmooth outputs due to the smoothness of their polynomial bases. The PCE-Kriging method from \cite{schobi2015polynomial} produces a more accurate and efficient approximation of highly nonlinear outputs compared to PCE. In PCE-Kriging, the PCE part approximates the global behavior of the output, while the Kriging part stochastically interpolates the local output behavior. However, computing the PCE component suffers from the curse of dimensionality, and constructing the PCE's basis functions is based on the assumption of independent inputs.    

Variance reduction methods (e.g., importance sampling \cite{srinivasan2002importance,peherstorfer2016multifidelity} or control variates~\cite{billinton1997composite}) exploit a problem-dependent sampling strategy to reduce the required number of expensive high-fidelity model evaluations. The efficiency of an importance-sampling method relies on the choice of importance-sampling density (or biasing density), which is often selected from a parameterized family of distributions. Importance-sampling methods can drastically reduce the number of samples required to estimate  risk~\cite{peherstorfer2016multifidelity} or reliability~\cite{tabandeh2022review}. To speed up the importance sampling process, the multifidelity importance sampling (MFIS) method from \cite{peherstorfer2016multifidelity} suggests using a surrogate model for learning the biasing density. In the context of MFIS for risk measure estimation, the authors in~\cite{HKTQ18CVaRROMS} used ROMs to construct a biasing density that samples in the \textit{risk region} of the parameter space, where input samples produce quantities of interest that are deemed ``risky''. The risk region is referred to as the $\epsilon$-risk region when it includes ROM error bounds. The work~\cite{HKTQ18CVaRROMS} shows that MFIS can estimate CVaR efficiently and without bias, even when the ROM has significant error compared to the high-fidelity model. 
However, that work does not consider high-dimensional or dependent random inputs, as the ROM construction becomes computationally intensive as the number of inputs increases.

This paper presents novel computational methods for CVaR estimation for nonlinear systems under high-dimensional and dependent inputs. First, we propose a novel \textit{DD-GPCE-Kriging} surrogate model, a fusion of DD-GPCE~\cite{lee2020practical} and Kriging to approximate a highly nonlinear and nonsmooth random output. Second, we leverage this new surrogate model in sampling-based CVaR estimation via (1) Monte Carlo sampling  and (2) multifidelity importance sampling. The proposed MCS-based method replaces the high-fidelity model (expensive to evaluate) with DD-GPCE-Kriging (cheap to evaluate) and is shown to be accurate in the presence of high-dimensional and dependent random inputs. However, the proposed method may produce a biased CVaR estimate as it relies on the approximation quality of the DD-GPCE-Kriging. Thus, we also suggest leveraging MFIS in connection with DD-GPCE-Kriging to determine the biasing density efficiently. The high-fidelity model is used for the importance-sampling-based CVaR estimate. 
For the proposed MFIS-based method, there are two major contributions. First, we suggest using the confidence interval (CI) or the probabilistic error bound of DD-GPCE-Kriging to define the $\epsilon$-risk region for the biasing density. Second, we compute DD-GPCE-Kriging by using computationally cheap low-fidelity output samples to further speed up the MFIS process. We investigate how much fidelity (and/or correlation with the high-fidelity model) is needed to construct a good biasing density.  
We demonstrate the performance of the two proposed CVaR estimation methods on lightweight composite structures with high stiffness and preferable heat resistance. It is essential to measure risk for these composites (used in high-performance air- and spacecraft) to produce risk-averse designs that can perform missions under harsh operating environments, e.g., high temperature and/or pressure. 
This work differs from our prior work~\cite{lee2022bi} as follows. First, we introduce the novel DD-GPCE-Kriging surrogate to improve the accuracy of the DD-GPCE for nonlinear outputs under dependent random inputs. Second, we suggest the MFIS-based method in connection with DD-GPCE-Kriging, which produces unbiased estimates and achieves higher efficiency than the MCS-based method in computing CVaR estimates.

The paper is organized as follows. Section~\ref{sec:2} covers the theoretical background, including the definition of the input and output random variables and introduces CVaR. We also summarize the DD-GPCE and necessary material on (multifidelity) importance sampling. Section~\ref{sec:3} presents the new DD-GPCE-Kriging surrogate model, and uses it for MCS-based CVaR estimation. Section~\ref{sec:4} proposes an MFIS-based method leveraging DD-GPCE-Kriging for efficient CVaR computation. Numerical results are reported in Section~\ref{sec:5}. In Section~\ref{sec:6}, we offer conclusions and an outlook for future work. 
\section{Theoretical background} \label{sec:2} 
This section presents the preliminaries for the $\mathrm{CVaR}$ problem setup in Section~\ref{sec:2.1}, defines input and output random variables in Sections~\ref{sec:2.2} and \ref{sec:2.3}, respectively, and then gives the $\rm{CVaR}$ definition in Section~\ref{sec:2.4}. We briefly summarize DD-GPCE in Section~\ref{sec:2.5} and review importance sampling in Section~\ref{sec:2.6} and multifidelity importance sampling in Section~\ref{sec:2.7}.
%
\subsection{Preliminaries} \label{sec:2.1} 
Let $\nat$, $\nat_{0}$, $\real$, and $\real_{0}^{+}$ be the sets of positive integers, non-negative integers, real numbers, and non-negative real numbers, respectively. For a positive integer $N\in\nat$, denote by $\mathbb{A}^N \subseteq \real^N$ a subspace of $\real^N$. 
\subsection{Input random variables} \label{sec:2.2} 
Let $(\Omega,\cF,\mathbb{P})$ be an abstract probability space, with a sample space $\Omega$, a $\sigma$-algebra $\cF$ on $\Omega$, and a probability measure $\mathbb{P}:\cF\to[0,1]$. Consider an $N$-dimensional random vector $\bX:=(X_{1},\ldots,X_{N})^\intercal:\Omega \rightarrow \mathbb{A}^N$, that models the uncertainties in a stochastic problem. We refer to $\bX$ as the random input vector or the input random variables. 
Denote by $F_{\bX}({\bx}):=\mathbb{P}\big[\cap_{i=1}^{N}\{ X_i \le x_i \}\big]$ the joint distribution function of $\bX$, admitting the joint probability density function $f_{\bX}({\bx}):={\partial^N F_{\bX}({\bx})}/{\partial x_1 \cdots \partial x_N}$. For $(\Omega,\cF,\mathbb{P})$, the image probability space is $(\mathbb{A}^N,\cB^{N},f_{\bX}(\bx)\rd\bx)$, where $\mathbb{A}^N$ is the image of $\Omega$ under the mapping $\bX:\Omega \to \mathbb{A}^N$ and $\cB^N:=\cB(\mathbb{A}^N)$ is the Borel $\sigma$-algebra on $\mathbb{A}^N\subset \mathbb{R}^N$.

\subsection{Output random variables}  \label{sec:2.3}
Given an input random vector $\bX$
with a known probability density function $f_{\bX}({\bx})$ on $\mathbb{A}^N \subseteq \real^N$, denote by $y(\bX)$ a real-valued, square-integrable  transformation on $(\Omega, \cF)$.  Here, $y:\mathbb{A}^N \to \real$ describes an output that an application engineer deems relevant for risk assessment. In this work, we assume that $y$ belongs to the weighted $L^2$ space $\left\{y:\mathbb{A}^N \to \real:~
		\int_{\mathbb{A}^N} \left| y(\bx)\right|^2  		f_{\bX}({\bx})\rd\bx < \infty \right\}$, which is the Hilbert space $\left\{Y=y(\bX):\Omega \to \real:
		~\int_{\Omega} \left|y(\bX(\omega))\right|^2 \rd\mathbb{P}(\omega) < \infty\right\}$ for the abstract probability space $(\Omega,\cF,\mathbb{P})$. If there is more than one output variable, then each component is associated with a measurement function $y_i$, $i=1,2,\ldots$.  The generalization for a multivariate output random vector is straightforward.
\subsection{Conditional Value-at-Risk (CVaR)} \label{sec:2.4} 
Given a random input vector $\bX$, consider an output function $y(\bX)\in L^2(\Omega,\cF,\mathbb{P})$. Let $\mathrm{CVaR}_{\beta}[y(\bX)]$ and $\mathrm{VaR}_{\beta}[y(\bX)]$ be the Conditional Value-at-Risk (CVaR) and the Value-at-Risk (VaR), respectively, of $y(\bX)$ at a given risk level $\beta\in(0,1)$. The $\mathrm{VaR}_{\beta}[y(\bX)]$ is the $\beta$-quantile of $y(\bX)$, i.e.,     
\begin{equation}
	\mathrm{VaR}_{\beta}[y(\bX)]=\argmin_{t\in\real}\{\mathbb{P}[y(\bX)\le t]\ge \beta\}.
	\label{2.2:1} 
\end{equation} 
Here $\mathbb{P}[y(\bX)\le t]=\int_{\mathbb{A}^N}\mathbb{I}_{\{y(\bx)\le t\}}(\bx)f_{\bX}(\bx)\rd\bx$, where the indicator function $\mathbb{I}_{\{y(\bx)\leq t\}}(\bx)$ is \emph{one}, if $y(\mathbf{x}\leq t)$, and \emph{zero}, otherwise. 
The $\mathrm{CVaR}_{\beta}[y(\bX)]$ is the mean of the outputs $y(\bX)$ exceeding $\mathrm{VaR}_{\beta}[y(\bX)]$.
At level $\beta\in(0,1)$, the $\mathrm{CVaR}_{\beta}$ can be determined as follows:
		\begin{align} 
			\mathrm{CVaR}_{\beta}[y(\bX)]=\mathrm{VaR}_{\beta}[y(\bX)]+\dfrac{1}{1-\beta}\Exp\left[(y(\bX)-\mathrm{VaR}_{\beta}[y(\bX)])_+\right],
			\label{2.2:3} 
		\end{align} 
		where $\Exp$ is the expectation operator with respect to $f_{\bX}(\bx)\rd\bx$, see~\cite*{rockafellar2000optimization,RTRockafellar_SUryasev_2002a}. If the cumulative distribution function $\mathbb{P}[Y\leq t]$ is continuous at $t=\mathrm{VaR}_{\beta}[y(\bX)]$, equation~\eqref{2.2:3} is equivalent to 
\begin{equation} 
	\mathrm{CVaR}_{\beta}[y(\bX)]=\dfrac{1}{1-\beta}\Exp[y(\bX)\cdot\mathbb{I}_{\cG_{\beta}[y(\bX)]}(\bX)],
	\label{2.2:4} 
\end{equation}  
where the \textit{risk region} for $\rm{CVaR}_{\beta}$ is defined as
\begin{align}\label{risk_region} 
	\cG_{\beta}[y(\bX)]:=\{\bx:y(\bx)\ge\mathrm{VaR}_{\beta}[y(\bX)]\}\subset\mathbb{A}^N. 	
\end{align}
With the definition of $\mathrm{CVaR}$ at hand, we can formalize the problem that we focus on in this work.

\begin{problem}\label{prob} 
	Consider a high-dimensional and dependent random input vector $\bX \in \mathbb{A}^N$. Given a highly nonlinear and/or nonsmooth output $y: \mathbb{A}^N\rightarrow \real$, requiring a computationally expensive model evaluation, our goal is to compute an unbiased estimate of $\mathrm{CVaR}_{\beta}[y(\bX)]$ accurately and efficiently.
\end{problem} 
\subsection{Dimensionally decomposed generalized polynomial chaos expansion}\label{sec:2.5}
For risk or reliability problems under high-dimensional inputs (e.g., $N\ge 20$), \emph{regular} GPCE~\cite{lee2020practical} or GPCE-based design methods \cite{lee2021rdo,lee2021rbdo} require a relatively large number of basis functions, thus suffering from the curse of dimensionality. Appendix~\ref{sec:appx1} provides details of \emph{regular} GPCE. 
This section briefly summarizes the more efficient DD-GPCE approximation~\cite{lee2023high}, which reorders and truncates the regular GPCE to  select a subset of the basis functions of the regular GPCE. The truncation is based on the degree of interaction among input variables, thereby tackling the curse of dimensionality to some extent. The chosen multivariate orthonormal polynomial basis functions are consistent with an arbitrary, non-product-type probability measure $f_{\bX}(\bx)\rd \bx$ of $\bX$. The bases can be determined by a three-step algorithm using a whitening transformation of the monomial basis, see Appendix~\ref{sec:appx2}.   

Given an input random vector $\bX$ with known probability density function $f_{\bX}(\bx)$, consider an output $y(\bX)$ as defined in Section \ref{sec:2.3}. 
To define the $S$-variate, $m$th-order DD-GPCE basis (see Appendix~\ref{sec:appx2} for details), we require $0 \le S \le N$ and $S \le m < \infty$ and define 
\begin{align} \label{2.4:1}  
	\bPsi_{S,m}(\bx):=(\Psi_1(\bx),\ldots,\Psi_{L_{N,S,m}}(\bx))^{\intercal}
\end{align} 
an $L_{N,S,m}$-dimensional column vector of orthonormal polynomials that are consistent with the probability measure $f_{\bX}(\bx)\rd \bx$. Here, the dimension of $\bPsi_{S,m}(\bx)$ is $L_{N,S,m}:=1+\sum_{s=1}^S\binom{N}{s}\binom{m}{s}$.
The $S$-variate, $m$th-order DD-GPCE approximation is then
\begin{equation}
	y_{S,m}(\bX) = \bc^{\intercal}\bPsi_{S,m}(\bX) \simeq y(\bX)
	\label{dd-gpce}
\end{equation}
with the column vector $\bc=(c_1,\ldots,c_{L_{N,S,m}})^{\intercal}$, whose elements are the expansion coefficients 
\begin{equation}
	\begin{split}
		c_i &:= \int_{\mathbb{A}^N} y(\bx)\Psi_i(\bx) f_{\bX}(\bx) \rd\bx,\quad i=1,\ldots,L_{N,S,m}.	\end{split}
\end{equation}
Based on \eqref{2.4:1}, a DD-GPCE approximation retains the degree of interaction among input variables less than or equal to $S$ and preserves polynomial orders less than or equal to $m$. 

\subsection{Importance sampling}\label{sec:2.6}
Importance sampling (IS) is a variance reduction method to improve the efficiency of MCS, and it follows two steps: (1) construct a biasing density, which is biased towards a specific region in input space (e.g., the risk region) and (2) sample from the biasing density and evaluate the high-fidelity model to estimate statistics of the output. The resulting IS estimator re-weights the drawn samples to compensate for the increased occurrence of biased samples in the sampling estimate~\cite{HKTQ18CVaRROMS}. 
To illustrate this for CVaR estimation, consider a random input vector $\bX$ with a probability density function $f_{\bX}(\bx)$ on $\mathbb{A}^N\subseteq\mathbb{R}^N$, consider the output $y(\bX)$. Let $f_{\bZ}(\bz)$ be a biasing density that defines a random vector $\bZ\in\bar{\mathbb{A}}^N$ such that $\mathbb{A}^N \subseteq \bar{\mathbb{A}}^N$. Consider the risk region \eqref{risk_region} corresponding to $\mathrm{CVaR}_{\beta}[y(\bX)]$, that is a small region in input space associated with the output exceeding $\mathrm{VaR}_{\beta}[y(\bX)]$. From \eqref{2.2:4}, we have  
		\begin{align}  
			\mathrm{CVaR}_{\beta}[y(\bX)]
			=\dfrac{1}{1-\beta}\mathbb{E}[\mathbb{I}_{\cG_{\beta}[y(\bZ)]}(\bZ)y(\bZ)w(\bZ)] \label{2.5:1:4}, 
		\end{align}  
		where $w(\bZ):=\frac{f_{\bX}(\bZ)}{f_{\bZ}(\bZ)}$ is a weight function and the expectation operator $\mathbb{E}$ is with respect to the biasing density $f_{\bZ}(\bz)$ of the random input vector $\bZ$. In IS, it is sufficient that the biasing density $f_{\bZ}(\bz) > 0$ for  $\bz\in\cG_{\beta}[y(\bX)]\subset \bar{\mathbb{A}}^N.$ The expectation in \eqref{2.5:1:4} can be approximated by MCS with samples $\bz_1,\ldots,\bz_M$, $M\in \mathbb{N}$, drawn from the biasing density $f_{\bZ}(\bz)$. Therefore, the IS estimator of $\mathrm{CVaR}_{\beta}[y(\bX)]$ is 
\begin{align}  \label{2.5:5}
	\widehat{\mathrm{CVaR}}_{\beta}^{\mathrm{IS}}[y(\bX)]=\dfrac{1}{M}\displaystyle\sum_{l=1}^M\mathbb{I}_{\cG_{\beta}[y(\bX)]}(\bz^{(l)})w(\bz^{(l)}).
\end{align} 
Compared to  MCS, one typically needs much fewer samples $M$ for IS to obtain converged $\mathrm{CVaR}$. 
In IS, a biasing density $f_{\bZ}(\bz)$ is designed so that it minimizes the variance of the IS estimator~\eqref{2.5:5}.  
\subsection{Multifidelity importance sampling}\label{sec:2.7}
Multifidelity importance sampling (MFIS) uses a surrogate model to design the biasing density while using the high-fidelity model to evaluate the IS estimator.
In~\cite{HKTQ18CVaRROMS}, the authors introduced a ROM-based MFIS method tailored to CVaR estimation. Given an error bound relating the output $y$ and its ROM approximation $\tilde{y}$, i.e, $|y(\bx)-\tilde{y}(\bx)|\leq\epsilon(\bx)$ for $\bx\in\mathbb{A}^N$, the $\epsilon$-\textit{risk region} is defined by
\begin{align} \label{2.6:1}
	\cG_{\beta}^{\epsilon}[\tilde{y}(\bX)]:=\{\bx:\tilde{y}(\bx)+\epsilon(\bx)\ge \mathrm{VaR}_{\beta}[\tilde{y}(\bX)-\epsilon(\bX)]\}.	
\end{align} 
By using an asymptotic result of the variance of $\widehat{\mathrm{CVaR}}_{\beta}^{\mathrm{IS}}[y(\bX)]$, \cite[Section~4]{HKTQ18CVaRROMS}, showed that the optimal biasing density and the weight function in the context of CVaR estimation are 
		\begin{align} 
			f_{\bZ}(\bz):=\dfrac{\mathbb{I}_{\cG_{\beta}^{\epsilon}[\tilde{y}(\bX)]}(\bz)f_{\bX}(\bz)}{\mathbb{P}[\cG_{\beta}^{\epsilon}[\tilde{y}(\bX)]]}\quad\text{and}\quad
			\label{2.6:2} 
			w(\bz):=\dfrac{\mathbb{P}[\cG_{\beta}^{\epsilon}[\tilde{y}(\bX)]]}{\mathbb{I}_{\cG_{\beta}^{\epsilon}[\tilde{y}(\bX)]}(\bz)}.
\end{align}
The MFIS estimator for CVaR then samples from $f_{\bZ}(\bz)$ and reweighs the samples according to \eqref{2.6:2}. 

\section{Dimensionally decomposed GPCE-Kriging for CVaR estimation} \label{sec:3} 
This section proposes a new DD-GPCE-Kriging method to efficiently and accurately estimate CVaR. Compared to the DD-GPCE introduced in Section \ref{sec:2.4}, the DD-GPCE-Kriging leverages both DD-GPCE \textit{and} Kriging. This allows for a more accurate approximation than DD-GPCE for highly nonlinear and nonsmooth outputs under dependent and high-dimensional random inputs.
A combination of standard PCE and Kriging was first proposed in \cite{schobi2015polynomial} in the context of reliability analysis. We propose DD-GPCE-Kriging in Section~\ref{sec:3.1} and give details on the calculation of the DD-GPCE-Kriging predictor in Section~\ref{sec:3.2}. In Section~\ref{sec:3.3}, we present an algorithm to estimate $\mathrm{CVaR}$ via DD-GPCE-Kriging.  

\subsection{DD-GPCE-Kriging model}\label{sec:3.1}
For $y(\bX)\in L^2(\Omega,\cal{F},\mathbb{P})$, we combine the $S$-variate, $m$th-order DD-GPCE approximation~\eqref{dd-gpce} with a Kriging model. The combination (referred to as DD-GPCE-Kriging) results in a Gaussian random variable of $y(\bx)$, i.e., 
\begin{equation} \label{DD-GPCE-Kriging} \bar{Y}_{S,m}(\bx)=\bc^{\intercal}\boldsymbol{\Psi}_{S,m}(\bx)+\sigma^2 Z(\bx;\btheta),
\end{equation}  
where the first and the second terms in the right-hand side of \eqref{DD-GPCE-Kriging} indicate the $S$-variate, $m$th-order DD-GPCE and Kriging, respectively. In \eqref{DD-GPCE-Kriging}, $\bc^{\intercal}\boldsymbol{\Psi}_{S,m}(\bx)$ follows the definition in \eqref{dd-gpce}, $\sigma^2$ is the Gaussian process variance, and $Z(\bx;\btheta)$ is a \emph{zero}-mean, \emph{unit}-variance stationary Gaussian process, which is fully determined by the autocorrelation function between two distinct input realizations $\bx$ and $\bx'$, i.e., 
\[
R(\bx,\bx';\btheta)=R(|\bx-\bx'|;\btheta).
\]
Here, $\btheta=(\theta_1,\ldots,\theta_N)^{\intercal}$ is the $N$-dimensional column vector of hyper-parameters to be computed. In this work, we chose either the Gaussian autocorrelation function 
\begin{equation}\label{auto-gauss} 
	R(|\bx-\bx'|;\btheta)=\exp\left[-\displaystyle\sum_{i=1}^{N}\left(\dfrac{x_i-x'_i}{\theta_i} \right)^2\right]
\end{equation}
or the exponential autocorrelation function 
\begin{equation}\label{auto-exp} 
	R(|\bx-\bx'|;\btheta)=\exp\left[-\displaystyle\sum_{i=1}^{N}\dfrac{|x_i-x'_i|}{\theta_i}\right].
\end{equation} 
We choose the autocorrelation function that most accurately predicts the quantity of interest for the particular application. There are also other choices (e.g., the Matern autocorrelation function), see~\cite{durrande2012additive,schobi2015polynomial}. In \eqref{auto-gauss} or \eqref{auto-exp}, an $N$-dimensional vector $\btheta$ of hyper-parameters can be determined by a maximum likelihood estimate~\cite{marrel2008efficient} or a leave-one-out cross-validation estimate~\cite{bachoc2013cross}. In this work, we select the leave-one-out cross-validation estimate method to evaluate all numerical examples since it provides more stable solutions than the maximum likelihood estimate method when the autocorrelation function family for outputs is unknown (see \cite{bachoc2013cross,schobi2015polynomial}). In Appendix \ref{sec:appx3}, we give details on the leave-one-out cross-validation estimate.   

\begin{remark} 
We can view DD-GPCE-Kriging as an extension of PCE-Kriging~\cite{schobi2015polynomial}. In contrast to the conventional PCE-Kriging, the DD-GPCE-Kriging method can handle arbitrary, dependent random inputs directly without relying on a detrimental transformation (e.g., Rosenblatt or Nataf transformation)
that may slow the convergence rate of the estimate to the true output. Moreover, the DD-GPCE part itself alleviates the curse of dimensionality to some extent by effectively truncating the basis functions and expansion coefficients in a dimensionwise manner.   
\end{remark} 

\subsection{Calculation of DD-GPCE-Kriging predictor}\label{sec:3.2}
Kriging is a stochastic interpolation method that provides the best predictor and an estimate of its accuracy. Given the known distribution of a random input $\bX$ and an output model $y:\mathbb{A}^N\rightarrow\mathbb{R}$, consider an input-output data set $\{\bx^{(l)},y(\bx^{(l)})\}_{l=1}^{L'}$ of size $L'\in\mathbb{N}$ constructed by evaluating the quantity of interest $y$ at each input data $\bx$, which are drawn as part of a sampling method (e.g., standard MCS, quasi MCS, or Latin hypercube sampling). We use the input-output data set to determine the $S$th-variate, $m$th-order DD-GPCE-Kriging at a new realization $\bx$. Then, as \textit{the best predictor} for $y(\bx)$, we have the mean $\bar{y}_{S,m}(\bx)$ of the Gaussian random variable $\bar{Y}_{S,m}$ in \eqref{DD-GPCE-Kriging}, i.e.,
\begin{equation} \label{mean}
	\bar{y}_{S,m}(\bx)=\hat{\bc}^{\intercal}\bPsi_{S,m}(\bx)+\br(\bx)^{\intercal}\bR^{-1}(\bb-\bA\hat{\bc}),
\end{equation} 
and the variance of $\bar{Y}_{S,m}$ is
\begin{align} \label{var} 
	\begin{split} 
		{\bar{\sigma}}_{S,m}^2(\bx)=&\hat{\sigma}^2\left\{1-\br(\bx)^{\intercal}\bR^{-1}\br(\bx)+\left[\bA^{\intercal}\bR^{-1}\br(\bx)-\bPsi_{S,m}(\bx)\right]^{\intercal}\left(\bA^{\intercal}\bR^{-1}\bA\right)^{-1} \right.\\
		&\left.\times\left[\bA^{\intercal}\bR^{-1}\br(\bx)-\bPsi_{S,m}(\bx)\right] \right\},
	\end{split} 
\end{align}
where
\begin{equation}
	\begin{split}
		\bA &:=
		\begin{bmatrix}
			\Psi_1(\bx^{(1)}) & \cdots &  \Psi_{L_{N,S,m}}(\bx^{(1)}) \\
			\vdots                   & \ddots &  \vdots                           \\
			\Psi_1(\bx^{(L^{\prime})}) & \cdots &  \Psi_{L_{N,S,m}}(\bx^{(L^{\prime})})
		\end{bmatrix}, \\
		\bR &:=
		\begin{bmatrix}
			R(\bx^{(1)},\bx^{(1)};\btheta) & \cdots &  R(\bx^{(1)},\bx^{(L')};\btheta) \\
			\vdots                   & \ddots &  \vdots                           \\
			R(\bx^{(L')},\bx^{(1)};\btheta) & \cdots &  R(\bx^{(L')},\bx^{(L')};\btheta) 
		\end{bmatrix},\\
	\end{split}
\end{equation} 
and $\br(\bx)=(R(\bx,\bx^{(1)};\btheta),\ldots,R(\bx,\bx^{(L')};\btheta))^{\intercal}$ and $\bPsi_{S,m}(\bx)=(\Psi_1(\bx),\ldots,\Psi_{L_{N,S,m}}(\bx))^{\intercal}$. In \eqref{mean}, $\bb=(y(\bx^{(1)}),\ldots,y(\bx^{(L')}))^{\intercal}$ is an $L'$-dimensional column vector of outputs evaluated at each input $\bx^{(l)}$, $l=1,\ldots,L'$, and $\hat{\bc}$ is an $L_{N,S,m}$-dimensional solution vector to the DD-GPCE's coefficients $\bc$ in \eqref{DD-GPCE-Kriging}, determined by
\begin{equation}\label{sol_c} 
	\hat{\bc} = \left(\bA^{\intercal}\bR^{-1}\bA\right)^{-1}\bA^{\intercal}\bR^{-1}\bb. 
\end{equation} 
In \eqref{var}, the variance term is 
\begin{equation}\label{sol_sig} 
	\hat{\sigma}^2 = \dfrac{1}{L'}\left(\bb-\bA\hat{\bc}\right)^{\intercal}\bR^{-1}\left(\bb-\bA\hat{\bc}\right).
\end{equation} 

In DD-GPCE-Kriging, the DD-GPCE model expands the output $y(\bx)$ into a set of polynomials to fit the mean behavior of the output. The Kriging model interpolates local variations with the weighted neighboring samples. As a result, DD-GPCE-Kriging leverages both DD-GPCE and Kriging, thus effectively predicting a highly nonlinear, nonsmooth output.       

\subsection{Sampling-based CVaR estimation by DD-GPCE-Kriging} \label{sec:3.3} 	
The estimation of $\mathrm{VaR}$~\eqref{2.2:1} and $\mathrm{CVaR}$~\eqref{2.2:3} for nontrivial examples requires a sampling method (e.g., standard MCS, quasi MCS, or Latin hypercube sampling). In this work, we follow \cite[Alg. 2.1]{HKTQ18CVaRROMS} and use standard MCS.
We use the predictor $\bar{y}_{S,m}$ of an $S$-variate, $m$th-order DD-GPCE-Kriging (see Section \ref{sec:3.2}) as a surrogate model to replace an expensive high-fidelity function $y(\bX)$. 
Therefore, we can afford a large number $L$ of output samples to obtain statistically converged CVaR estimates when using the DD-GPCE-Kriging surrogate. 
The sampling-based estimation method evaluates $\bar{y}_{S,m}(\bx)$, and we denote the resulting estimates by $\widehat{\mathrm{VaR}}_{\beta}[\bar{y}_{S,m}(\bX)]$ and $\widehat{\mathrm{CVaR}}_{\beta}[\bar{y}_{S,m}(\bX)]$.

Algorithm~\ref{al1} summarizes the steps of the proposed $\mathrm{CVaR}$ estimation process.
Estimating $\mathrm{VaR}$ is straightforward (see Steps 3-4).
From the definition of ${\mathrm{CVaR}}_{\beta}[y(\bX)]$ in \eqref{2.2:3}, we obtain the estimate 
\begin{align} 
	\widehat{\mathrm{CVaR}}_{\beta}[\bar{y}_{S,m}(\bX)]=\widehat{\mathrm{VaR}}_{\beta}[\bar{y}_{S,m}(\bX)]+\dfrac{1}{1-\beta}\dfrac{1}{L}\displaystyle\sum_{l=1}^L\left(\bar{y}_{S,m}(\bx^{(l)})-\widehat{\mathrm{VaR}}_{\beta}[\bar{y}_{S,m}(\bX)]  \right)_+.
	\label{3.2}	
\end{align} 
We use \eqref{3.2} in the final step of Algorithm~\ref{al1}.

\begin{algorithm}
	\caption{Sampling-based estimation of $\mathrm{CVaR}_{\beta}$ by the DD-GPCE-Kriging predictor $\bar{y}_{S,m}$.}
	\label{al1}
	\begin{algorithmic}[1] 
		\Require Input samples $\bx^{(l)}=(x_1^{(l)},\ldots,x_N^{(l)})^{\intercal}$, $l=1,\ldots, L,~L\gg1$, via MCS, quasi MCS, or Latin hypercube sampling with corresponding probabilities $p^{(l)}=f_{\bX}(\bx^{(l)}){\rm d} \bx^{(l)}$;  Set sample number $L'<L$, used for the computation of DD-GPCE-Kriging; Set a risk level $\beta \in (0,1)$.
		\Ensure  Estimates $\widehat{\mathrm{VaR}}_{\beta}[\bar{y}_{S,m}(\bX)]$ and $\widehat{\mathrm{CVaR}}_{\beta}[\bar{y}_{S,m}(\bX)]$.
		\Procedure{Calculate DD-GPCE-Kriging}{$\bar{y}_{S,m}(\bx)$}
		\State Create measure-consistent orthonormal basis vector $\bPsi_{S,m}(\bx)$. We refer to Appendix~\ref{sec:appx2}.
		\State Determine hyper-parameters $\btheta$ by either of a maximum likelihood estimate or a leave-one-out cross-validation estimate, see Section~\ref{sec:3.1} and Appendix~\ref{sec:appx3}.
		\State Calculate DD-GPCE's coefficients $\bc$~\eqref{sol_c}.
		\EndProcedure 
		\Procedure{Create output samples}{$\{\bar{y}_{S,m}(\bx^{(l)})\}_{l=1}^{L}$}
		\State Create $L$ output samples of the predictor $\bar{y}_{S,m}$~\eqref{mean} of DD-GPCE-Kriging for $y$.
		\EndProcedure 
		\State Sort values of $\bar{y}_{S,m}$ in descending order and relabel the samples so that
		$$
		\bar{y}_{S,m}(\bx^{(1)})>\bar{y}_{S,m}(\bx^{(2)})>\cdots>\bar{y}_{S,m}(\bx^{(L)}),$$
		and reorder the probabilities accordingly (so that $p^{(l)}$ corresponding to $\bx^{(l)}$).
		\State Compute the index $\tk_{\beta}\in\nat$ such that 
		$$
		\displaystyle\sum_{l=1}^{\tk_{\beta}-1}p^{(l)} \le 1-\beta < \sum_{l=1}^{\tk_{\beta}}p^{(l)}. 
		$$
		\State Set $\mathrm{VaR}_{\beta}[y(\bX))]\approx\widehat{\mathrm{VaR}}_{\beta}[\bar{y}_{S,m}(\bX)]={\bar{y}_{S,m}}(\bx^{({\tilde{k}}_{\beta})})$.
		\State Set $\mathrm{CVaR}_{\beta}[y(\bX)]\approx \widehat{\mathrm{CVaR}}_{\beta}[\bar{y}_{S,m}(\bX)]$~\eqref{3.2}.  
	\end{algorithmic}
\end{algorithm}

To satisfy the accuracy of the DD-GPCE-Kriging's predictor, we need to determine the best-fit parameters (e.g., $\hat{\bc}$, $\hat{\sigma}$, or $\btheta$ in \eqref{sol_c} and \eqref{sol_sig}). We select the number of output samples to be at least greater than the number $L_{N,S,m}$ of DD-GPCE's coefficients or basis functions. Consequently, when faced with high-dimensional inputs (e.g., $N\ge20$), DD-GPCE-Kriging may require thousands of output samples. In such cases, computing a DD-GPCE-Kriging surrogate may require prohibitive computational cost when each sample is determined by an expensive high-fidelity model evaluation. Moreover, when not using enough samples in the computation of DD-GPCE-Kriging, a significant bias in the CVaR estimate may be present. In the following section, we propose to use the DD-GPCE-Kriging surrogate as a part of an MFIS framework. This further reduces the computational expense, while obtaining unbiased CVaR estimates via importance sampling.

\section{Multifidelity importance sampling for CVaR estimation} \label{sec:4}  
We present a novel multifidelity importance sampling strategy that leverages the DD-GPCE-Kriging surrogates (tailored to CVaR) to obtain an unbiased CVaR estimate for a highly nonlinear and nonsmooth random output. In contrast to existing methods~\cite{peherstorfer2016multifidelity,HKTQ18CVaRROMS,lee2022bi}, the proposed MFIS-based method works for high-dimensional and dependent random inputs (see Problem~\ref{prob}). 
In the proposed MFIS-based method, we use the DD-GPCE-Kriging surrogate to determine the biasing density, from which we then draw a few ($M \ll L$) high-fidelity output samples to estimate $\mathrm{VaR}_{\beta}$ and $\mathrm{CVaR}_{\beta}$. We use the confidence interval (CI) of the predictor $\bar{y}_{S,m}$ to estimate the risk region for an optimal biasing density~\eqref{2.6:2} which is determined via the CI-based  $\epsilon$-risk region. 
We define the CI-based $\epsilon$-risk region in Section~\ref{sec:4.1}, present the importance sampling in Section~\ref{sec:4.2}, and introduce the complete algorithm to obtain unbiased $\mathrm{CVaR}_{\beta}$ estimates in Section~\ref{sec:4.3}.  
\subsection{Confidence interval based $\epsilon$-risk region} \label{sec:4.1}
\begin{algorithm}
	\caption{CI-based $\epsilon$-risk region by the DD-GPCE-Kriging predictor $\bar{y}_{S,m}$.}
	\label{al2}
		\begin{algorithmic}[1] 
			\Require Samples $\bx^{(l)}=(x_1^{(l)},\ldots,x_N^{(l)})^{\intercal}$, $l=1,\ldots, L$, via MCS, quasi MCS, or Latin hypercube sampling with corresponding probabilities $p^{(l)}=f_{\bX}(\bx^{(l)}){\rm d} \bx^{(l)}$; Generate the input-output data set $\{\bx^{(l)},y(\bx^{(l)})\}_{l=1}^{\bar{L}}$, $\bar{L}<L$; Set risk level $\beta\in(0,1)$.
			\Ensure  CI-based $\epsilon$-risk region $\widehat{\cG}_{\beta}^{\alpha}$.
			\State Calculate the $S$-variate, $m$th-order DD-GPCE-Kriging $\bar{Y}_{S,m}$ through Steps 1--5 of Algorithm \ref{al1}. 
			\State  Evaluate the mean value $\bar{y}_{S,m}(\bx^{(l)})$ \eqref{mean} and the standard deviation $\bar{\sigma}_{S,m}(\bx^{(l)})$ \eqref{var} of $\bar{Y}_{S,m}$ at samples $\{\bx^{(l)}\}_{l=1}^{L}$ and determine the corresponding the lower or upper limits of CI, $\epsilon_{S,m}^{\alpha}(\bx^{(l)})$~\eqref{4.1:2}. 
			\State Sort values of $(\bar{y}_{S,m}(\bx^{(l)})-\epsilon_{S,m}({\bx}^{(l)}))$ in descending order and relabel those, i.e.,
			\[
			\bar{y}_{S,m}(\bx^{(1)})-\epsilon_{S,m}^{\alpha}(\bx^{(1)})>\cdots>\bar{y}_{S,m}(\bx^{(L)})-\epsilon_{S,m}^{\alpha}(\bx^{(L)}).
			\] 
			\State Compute an index $\bar{k}_{\beta}$ such that 
			\[
			\sum_{l=1}^{\bar{k}_{\beta}-1}p^{(l)}\leq 1-\beta < \sum_{l=1}^{\bar{k}_{\beta}}p^{(l)}.
			\]
			\State Set
			\[
			\widehat{\mathrm{VaR}}_{\beta}[\bar{y}_{S,m}(\bX)-\epsilon_{S,m}^{\alpha}(\bX)]=\bar{y}_{S,m}(\bx^{(\bar{k}_{\beta})}).
			\]
			\State Determine the $(1-\alpha)100\%$ CI-based $\epsilon$-risk region by a discrete set: 
			\[
			\widehat{\cG}_{\beta}^{\alpha}:=\{\bx^{(j)} \in \mathbb{A}^N~:~\bar{y}_{S,m}(\bx^{(j)})+\epsilon_{S,m}^{\alpha}(\bx^{(j)})\geq \widehat{\mathrm{VaR}}_{\beta}[\bar{y}_{S,m}(\bX)-\epsilon_{S,m}^{\alpha}(\bX)],~j=1,\ldots ,\bar{L}<L\}.
			\] 
		\end{algorithmic}
	\end{algorithm}
	The random variable $\bar{Y}_{S,m}$ of DD-GPCE-Kriging is assumed to follow the Gaussian distribution with the mean $\bar{y}_{S,m}(\bx)$ in \eqref{mean} and variance $\bar{\sigma}_{S,m}^2(\bx)$ in \eqref{var}. Given $\alpha\in[0,1]$, the upper and lower limits of the $(1-\alpha)100\%$ confidence interval of the DD-GPCE-Kriging predictor $\bar{y}_{S,m}$ at input $\bx$ satisfy 
	\begin{equation}
		-Q_{1-\alpha/2}{\bar{\sigma}}_{S,m}(\bx) \leq y(\bx)-\bar{y}_{S,m}(\bx) \leq Q_{1-\alpha/2}{\bar{\sigma}}_{S,m}(\bx).
		\label{CIboundary} 
	\end{equation}  
	Here, $Q_{1-\alpha/2}$ is the $(1-\alpha/2)$-level quantile of the standard normal distribution and ${\bar{\sigma}_{S,m}}^2(\bx)$ is the predicted variance at $\bx$ by DD-GPCE-Kriging. By expressing $Q_{1-\alpha/2}\bar{\sigma}_{S,m}(\bx)$ as ${\epsilon}_{S,m}^{\alpha}(\bx)$, equation~\eqref{CIboundary} becomes
	\begin{equation}\label{4.1:2} 
		-{\epsilon}_{S,m}^{\alpha}(\bx)\leq y(\bx)-\bar{y}_{S,m}(\bx)\leq {\epsilon}_{S,m}^{\alpha}(\bx).
	\end{equation}  
	
	In reference to Section \ref{sec:2.6}, by replacing $\epsilon$ in \eqref{2.6:1} with $\epsilon_{S,m}^{\alpha}$ in \eqref{4.1:2}, the $\epsilon$-risk region corresponding to $\mathrm{CVaR}_{\beta}[y(\bX)]$ becomes the CI-based $\epsilon$-risk region
	\begin{align}\label{erisk}
		\cG_{\beta}^{\alpha}:=\{\bx \in \mathbb{A}^N~:~\bar{y}_{S,m}(\bx)+\epsilon_{S,m}^{\alpha}(\bx)\geq\mathrm{VaR}_{\beta}[\bar{y}_{S,m}(\bX)-\epsilon_{S,m}^{\alpha}(\bX)]\}.
	\end{align}
	The procedure for computing the CI-based $\epsilon$-risk region~\eqref{erisk} is presented in Algorithm~\ref{al2}. 
	
	\begin{remark}
		In contrast to the previous work~\cite{HKTQ18CVaRROMS}, the CI (as a probabilistic bound) does not guarantee that the true risk region is included in the CI-based $\epsilon$-risk region at all times. 
		We note that the $\epsilon$-risk region is only used to define the biasing density. Thus, since the MFIS estimator uses high-fidelity output samples, it is always unbiased, even if the CI-based $\epsilon$-risk region is inaccurate. The quality of the $\epsilon$-risk region only impacts the variance of the CVaR estimate.  
	\end{remark}
	\subsection{Importance sampling for CVaR estimation} \label{sec:4.2}
	For importance sampling, we can determine an optimal biasing density by minimizing the asymptotic variance of the $\mathrm{CVaR}_{\beta}$ estimator (e.g., see Section 4.2 of \cite{HKTQ18CVaRROMS}). We use the CI-based $\epsilon$-risk region to construct the optimal biasing density as 
	\begin{align}\label{4.2:1}
		f_{\bZ}(\bz) := \dfrac{\mathbb{I}_{\cG_{\beta}^{\alpha}[\bar{y}_{S,m}(\bX)]}(\bz)f_{\bX}(\bz)}{\mathbb{P}[\cG_{\beta}^{\alpha}[\bar{y}_{S,m}(\bX)]]},
	\end{align}
	where $\mathbb{I}_{\cG_{\beta}^{\alpha}[\bar{y}_{S,m}(\bX)]}(\bz)$ is the indicator function associated with the CI-based $\epsilon$-risk region. 
	Inserting $f_{\bZ}(\bz)$ from equation~\eqref{4.2:1} into the weight function \eqref{2.6:2} leads to 
	\begin{align}\label{4.2:2}
		w(\bz) := \dfrac{\mathbb{P}[\cG_{\beta}^{\alpha}[\bar{y}_{S,m}(\bX)]]}{\mathbb{I}_{\cG_{\beta}^{\alpha}[\bar{y}_{S,m}(\bx)]}(\bz)}.
	\end{align}
	We draw a few high-fidelity output samples $\{y(\bz^{(l)})\}_{l=1}^{M}$, where $M\ll L$, from \eqref{4.2:1} and use \eqref{4.2:2} and \eqref{2.5:5} to obtain the importance sampling estimate of $\mathrm{CVaR}_{\beta}$, namely $\widehat{\mathrm{CVaR}_{\beta}}^{\rm{IS}}[y(\bX)]$~\eqref{2.5:5}.

	\subsection{Proposed multifidelity importance sampling for CVaR estimation} \label{sec:4.3}

	Algorithm~\ref{al3} outlines the procedure of the proposed MFIS-based method to obtain $\widehat{\mathrm{CVaR}}_{\beta}^{\mathrm{IS}}[y(\bX)]$, including Algorithms \ref{al2} and \ref{al3}. In Step~1 of Algorithm~\ref{al3}, we suggest computing the DD-GPCE-Kriging surrogate with two options:
	\begin{itemize} 
		\item[(a)] HF (High-fidelity) that uses the high-fidelity model evaluations, 
		\item[(b)] LF (Low-fidelity) that uses the low-fidelity model evaluations. 
	\end{itemize}

	We include here the low-fidelity model to speed up the computation of DD-GPCE-Kriging. We investigate how much fidelity/correlation of a surrogate model and the high-fidelity model is really needed to obtain good biasing densities. Oftentimes, the only option one may have is to use a coarse-grid approximation in order to generate enough data to train a surrogate model. We emphasize again that the proposed MFIS-based method with LF can provide an unbiased CVaR estimate even when the low-fidelity model has a significant error compared to the high-fidelity model. Absolute error, as it turns out, is not the sole driving factor, as \textit{correlation} with the high-fidelity output is a key driver, see also \cite{peherstorfer2016optimal}. We demonstrate this with a variety of (artificially created) low-fidelity models on a mathematical example in Section~\ref{sec:5}.
	
	In Algorithm~\ref{al3}, the cost to determine the CVaR estimate is dominated by constructing two sets of input-output data of size $L'$ and $M$. For the first input-output data set, its size $L'$ is determined as a multiple of the number of DD-GPCE's coefficients or basis functions $L_{N,S,m}$, satisfying the necessary condition ($L'\leq L_{N,S,m}$). In the proposed algorithm of MFIS, the computational cost to determine the first input-output data set can be reduced by replacing a high-fidelity output model with its low-fidelity version. For the second input-output data set, its sample number $M$ is usually chosen to be less than $L'$, and in this work, the chosen number is $L'/2$, but one can select a different number via a convergence test for the CVaR estimator in other problems.

	Let $c_T\in\mathbb{R}_0^+$ be the total computational budget and $c_H\in\mathbb{R}_0^+$ and $c_L\in\mathbb{R}_0^+$ be the costs of the high-fidelity and low-fidelity model evaluations, respectively. Then, the total cost $c_T$ of Algorithm~\ref{al3} is, in the HF option
	\[
	c_T = (L' + M)c_H 
	\]  
	and, in the LF option
	\[
	c_T = L'c_L + Mc_H.
	\]
	Given a total computational budget $c_T$ for LF, we can determine the cost of a low-fidelity model evaluation such that
	\[
	c_L\leq\dfrac{c_T-Mc_H}{L'}.
	\]  
	
	Generally, $M < L'$ and $c_L < c_H$, therefore, the proposed MFIS-based method with LF has higher efficiency than its HF version. More optimal use of the computational budget can be achieved, see~\cite{yang2022control}. 
	
	\begin{algorithm}
		\caption{Proposed MFIS-based method for $\mathrm{CVaR}_{\beta}$ estimation.}
		\label{al3}
		\begin{algorithmic}[1] 
			\Require Samples $\bx^{(l)}=(x_1^{(l)},\ldots,x_N^{(l)})^{\intercal}$, $l=1,\ldots,L,~L\gg1$, via MCS, quasi MCS, or Latin hypercube sampling with corresponding probabilities $p^{(l)}=f_{\bX}(\bx^{(l)}){\rm d}{\bx}^{(l)}$; Set a risk level $\beta \in (0,1)$.
			\Ensure  $\widehat{\mathrm{CVaR}}_{\beta}^{\mathrm{IS}}[\bar{y}_{S,m}(\bX)]$.
			\State Calculate the $S$-variate, $m$th-order DD-GPCE-Kriging $\bar{Y}_{S,m}$ through Steps 1--5 of Algorithm~\ref{al1}.   
			For Steps 3 and 4 of Algorithm~\ref{al1}, 
			\begin{itemize}
				\item[(a)] HF (High-fidelity): Use high-fidelity outputs of size $L'$ (e.g., fine mesh model for FEA).
				\item[(b)] LF (Low-fidelity): Use low-fidelity outputs of size $L'$ (e.g., coarse mesh model for FEA).
			\end{itemize}
			\State Determine the $(1-\alpha)100\%$ CI-based $\epsilon$-risk region by a discrete set through Steps 2--6 of Algorithm~\ref{al2}.    
			\State Select input samples $\bx^{(l)}$, $l=1,\ldots,M$, within the CI-based $\epsilon$-risk region. Then, construct high-fidelity output samples of $y$ at $\bx^{(l)}$.
			\State Assign values of the weight function \eqref{4.2:2}, i.e., $w(\bx^{(1)}),\ldots,w(\bx^{(M)})\equiv\mathbb{P}[\widehat{\cG}_{\beta}^{\alpha}[y(\bX)]]\approx \big|\widehat{\cG}_{\beta}^{\alpha}[y(\bX)]\big|/L$, where $\big|\cdot\big|$ denotes the cardinality of a set. 
			\State Obtain  $\mathrm{CVaR}_{\beta}[y(\bX)]\approx \widehat{\mathrm{CVaR}}_{\beta}^{\mathrm{IS}}[\bar{y}_{S,m}(\bX)]$ from \eqref{2.5:5} with $y(\bx^{(l)})$, $p^{(l)}=w(\bx^{(l)})/M$, $l=1,\ldots,M$.  
		\end{algorithmic}
	\end{algorithm}
	
\begin{figure*}
		\begin{center}
			\includegraphics[angle=0,scale=1.0,clip]{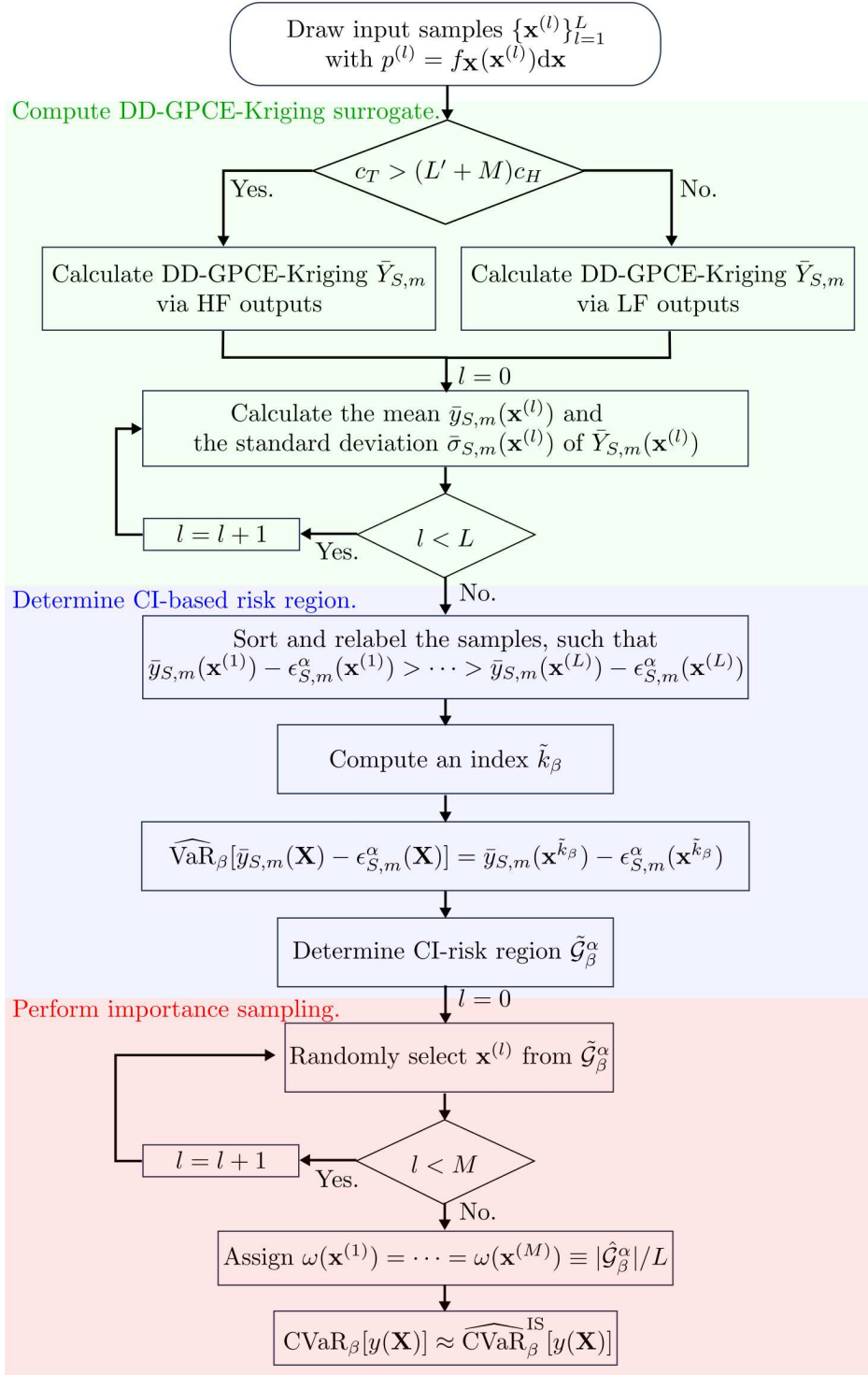}
		\end{center}
		\caption{Flow chart for the proposed MFIS-based method to estimate Conditional Value-at-Risk.}
		\label{fig1}
	\end{figure*}

	Figure~\ref{fig1} presents a flow chart for the proposed MFIS-based method. It starts with drawing the samples for the random input vector $\bX$. Given the computational budget $c_T$, we choose to evaluate either the high-fidelity or low-fidelity model to calculate DD-GPCE-Kriging according to $c_T \leq (L' + M)c_H$. We calculate the mean $\bar{y}_{S,m}(\bx^{(l)})$ and the standard deviation $\bar{\sigma}_{S,m}(\bx^{(l)})$ of $\bar{Y}_{S,m}(\bx^{(l)})$ for $l=1,\ldots,L$. From these, we estimate $\widehat{\rm{VaR}_{\beta}}[\bar{y}_{S,m}(\bX)-\epsilon_{S,m}^{\alpha}(\bX)]$ and then determine the CI-risk region $\tilde{G}_{\beta}^{\alpha}$. We conduct importance sampling by selecting $M$ input samples from $\tilde{G}_{\beta}^{\alpha}$. The importance sampling process produces the CVaR estimate $\widehat{\mathrm{CVaR}}_{\beta}^{\mathrm{IS}}[y(\mathbf{X})]$ with only $M$ high-fidelity outputs.
	
	The next section demonstrates the computational efficiency of the proposed MCS and MFIS-based methods via practical engineering examples. 
\section{Numerical results} \label{sec:5} 

Four numerical examples are presented to evaluate the proposed DD-GPCE-Kriging-based MCS and MFIS methods for $\mathrm{CVaR}$ estimation. Sections~\ref{sec:5.1} and \ref{sec:5.2} describe the details of the numerical setup and error measures, respectively, used for Examples~1--4. In Sections~\ref{sec:5.3} and \ref{sec:5.4}, our instructional examples consider two closed-form mathematical functions for which we construct a suite of surrogate models analytically.   
In Sections~\ref{sec:5.5} and \ref{sec:5.6}, we demonstrate the scalability of the proposed methods and their applicability to complex engineering problems via a 2D glass/vinylester composite plate example and a 3D composite T-joint example, respectively.  
\subsection{Numerical setup} \label{sec:5.1} 

In Examples~1--3, the sample number for MCS is $L=10,000$ and in Example~4, the sample number is $L=6,000$ due to the higher computational cost associated with that model. In all four examples, the chosen sample numbers led to converged CVaR estimates. We deem estimates converged if the difference between the previous and current CVaR value differs by less than 0.1\%. 

In all four examples, the monomial moment matrix $\bG_{S,m}$ (see \eqref{appx2:2} in Appendix~\ref{sec:appx2}) used for creating $\bPsi_{S,m}(\bx)$ in DD-GPCE-Kriging~\eqref{DD-GPCE-Kriging} is determined by quasi MCS with $5\times 10^6$ samples together with the Sobol sequence. To determine the hyper-parameters $\btheta$ via a leave-one-out cross-validation (presented in Appendix~\ref{sec:appx3}), in all examples, we employed the trust-region reflective optimization algorithm, which is known to perform well for nonlinear optimization problems~\cite{Conn2009opt}. We obtained all numerical results via MATLAB~\cite{Matlab21} on an Intel Core i7-10850H 2.70 GHz processor with 64 GB of RAM.  

\subsection{Error measures} \label{sec:5.2} 
To measure the deviation of a CVaR estimate relative to its standard MCS estimate (which we consider the benchmark), we use the mean relative difference (MRD). The MRD is an average error measure from $K\in\mathbb{N}$ independent estimates, i.e.,  
\begin{align}\label{mrd}
	\mathrm{MRD}=&\dfrac{\dfrac{1}{K}\displaystyle\sum_{k=1}^K\bigg|{\cal{Y}}_k-\widehat{\mathrm{CVaR}}_{\beta}[y(\bX)]   \bigg|}{\bigg|\widehat{\mathrm{CVaR}}_{\beta}[y(\bX)] \bigg|},
\end{align} 
where, ${\cal{Y}}_k=\widehat{\rm{CVaR}}_{\beta}^{(k)}[\bar{y}_{S,m}(\bX)]$ or $\widehat{\rm{CVaR}}_{\beta}^{\rm{IS}(k)}[y(\bX)]$ is the estimate at the $k$th independent trial run for DD-GPCE-Kriging-based MCS or MFIS, respectively. We also provide the normalized root-mean-square deviation (N-RMSD) as   
\begin{align}\label{msd}
\mathrm{N\text{-}RMSD}=\sqrt{\dfrac{\dfrac{1}{K}\displaystyle\sum_{k=1}^K\bigg({\cal{Y}}_k-\widehat{\mathrm{CVaR}}_{\beta}[y(\bX)]   \bigg)^2}{\left(\widehat{\mathrm{CVaR}}_{\beta}[y(\bX)]\right)^2}},
\end{align}
which measures the standard deviation of a CVaR estimator relative to the benchmark. The N-RMSD has a normalized value to the benchmark.

Additionally, in Examples 1, 3, and 4,  we determine the Pearson correlation coefficient (PCC) between the high-fidelity model output $y_H$ and low-fidelity model output $y_L$ from $P\gg 1$ samples, i.e.,
\begin{align}\label{pcor}
	\rho_{H,L}=\dfrac{\dfrac{1}{P}\displaystyle\sum_{l=1}^P\left(y_H(\bx^{(l)})-\dfrac{1}{P}\sum_{l=1}^P y_H(\bx^{(l)})\right)\left(y_L(\bx^{(l)})-\dfrac{1}{P}\sum_{l=1}^P y_L(\bx^{(l)})\right)}{\sqrt{\dfrac{1}{P}\displaystyle\sum_{l=1}^P \left(y_H(\bx^{(l)})-\dfrac{1}{P}\sum_{l=1}^P y_H(\bx^{(l)})\right)^2{\dfrac{1}{P}\displaystyle\sum_{l=1}^P\left(y_L(\bx^{(l)})-\dfrac{1}{P}\sum_{l=1}^P y_L(\bx^{(l)})\right)^2}}},
\end{align} 
where $\rho_{H,L}\in[-1,1]$. In Examples 1, 3, and 4, the selected values of $P$ in \eqref{pcor} are $10,000$, $3,000$, and $3,000$, respectively. The correlation coefficient $\rho_{H,L}$ guides us in selecting a low-fidelity model $y_L$ for computing DD-GPCE-Kriging in the proposed MFIS-based method with LF, similar to \cite{peherstorfer2016optimal}. We ignore a low-fidelity model when its correlation coefficient with the high-fidelity model is close to \emph{zero} (i.e., $|{\rho}_{H,L}|\ll 1$). In Example~1, we provide four distinct low-fidelity models $y_{L,i}$ and we compute the corresponding $\rho_{H,L_i}$ by replacing $y_L$ with $y_{L_i}$ in \eqref{pcor} for $i=1,2,3,4$. The following example illustrates how retaining qualitative features of the high-fidelity model in the low-fidelity model can impact the accuracy of the estimator $\widehat{\rm{CVaR}}_{\beta}^{\rm {IS}}$.

\subsection{Example 1: Rastrigin function and its surrogates} \label{sec:5.3}
This first example illustrates the DD-GPCE-Kriging-based MCS and MFIS methods for a highly nonlinear analytical mathematical function for which we construct simple surrogates. 

\begin{figure*}
	\begin{center}
		\includegraphics[angle=0,scale=1.0,clip]{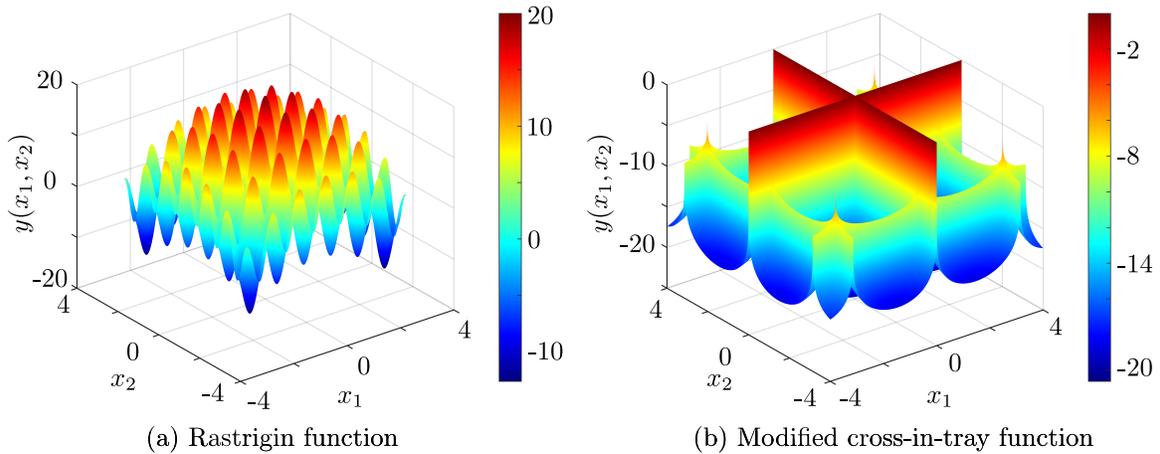}
	\end{center}
	\caption{Graphs of two two-dimensional functions evaluated for $\mathrm{CVaR}_{\beta}$ in Examples 1 and 2. The Rastrigin function in Figure~\ref{fig2}a exhibits high nonlinearity yet is smooth, while the modified cross-in-tray function in Figure~\ref{fig2}b exhibits high nonlinearity and nonsmoothness. Each color map indicates the values of $y(x_1,x_2)$ on the domain $[-4,4]\times[-4,4]$.}
	\label{fig2}
\end{figure*}

\subsubsection{Problem definition} 

For a bivariate Gaussian random vector $\bX=(X_1,X_2)^{\intercal}$ with mean vector $(0,0)^{\intercal}$ and standard deviations $\sqrt{\mathrm{var}[X_i]}=2,~i=1,2$, consider the Rastrigin function

\begin{align} \label{ras_f}
	y(X_1,X_2)=10-\sum_{i=1}^2(X_i^2-5\cos(2\pi X_i)).
\end{align} 
We consider two distinct cases of $\bX$, where the correlation coefficient between $X_1$ and $X_2$ is $0.9$ in the first case, and zero in the second case. Since the random input variables follow a Gaussian distribution, the \emph{zero} correlation further means that they are independent.

Figure~\ref{fig2}a shows the Rastrigin function on the domain $[-4,4]\times[-4,4]$. We see that the Rastrigin function has smooth, continuous, and highly oscillatory behavior. The function is notoriously difficult to approximate with a polynomial basis as it requires high-degree polynomial bases.

\subsubsection{High-fidelity and low-fidelity models}
This example seeks to illustrate that correlation between the low-fidelity and high-fidelity model is instrumental to construct good biasing densities. We consider the Rastrigin function as the high-fidelity quantity of interest (`truth'). We propose the following four low-fidelity models:    
\begin{itemize} 
	\item[(1)] Low-fidelity model \#1 (Output offset by a constant 90): 
	\begin{align}\label{sur1}
		y_{L_1}(\bX)=100-\sum_{i=1}^2(X_i^2-5\cos(2\pi X_i))=y(\bX)+90,
	\end{align}
	\item[(2)] Low-fidelity model \#2 (Magnification of the output by a factor of 10):
	\begin{align}\label{sur2}
		y_{L_2}(\bX)=100-\sum_{i=1}^2(10X_i^2-50\cos(2\pi X_i))=10y(\bX),
	\end{align}
	\item[(3)] Low-fidelity model \#3 (Phases shifted by $\pi/2$):
	\begin{align}\label{sur3}
		y_{L_3}(\bX)=10-\sum_{i=1}^2(X_i^2-5\cos(2\pi X_i+\pi/2)),
	\end{align}
	\item[(4)] Low-fidelity model \#4 (Frequencies scaled by 1/2):
	\begin{align}\label{sur4}
		y_{L_4}(\bX)=10-\sum_{i=1}^2(X_i^2-5\cos(\pi X_i)).
	\end{align}	 
\end{itemize}  
The correlation coefficients $\rho_{H,{L_i}}$ from \eqref{pcor} for $i=1,2,3,4$ between the high-fidelity model output $y_H$ and the low-fidelity model output $y_{L_i}$ are 1, 1, 0.72, and 0.72, respectively. These indicate that the low-fidelity models 1 and 2 have a stronger relationship with the high-fidelity model than the low-fidelity models 3 and 4. A perfect correlation of \emph{one} implies that the low-fidelity models 1 and 2 capture the trends of the high-fidelity model output completely, even when the low-fidelity models 1 and 2 have significant bias. We use the above four surrogates for MFIS to study how correlation and accuracy of the low-fidelity models relative to the high-fidelity model affect the estimator. 

\subsubsection{Results}

\begin{table*}
	\caption{$\mathrm{CVaR}_{\beta}$ estimates ($\beta=0.99$) for the Rastrigin function~\eqref{ras_f} via the DD-GPCE-Kriging, DD-GPCE, and PCE-Kriging-based MCS methods, the DD-GPCE-Kriging-based MFIS, and standard MCS.}
	\vspace{-0.1in}
	\begin{spacing}{1.2}
		\begin{center}
			\small
			\resizebox{\textwidth}{!}{
				\begin{tabular}{clccccccc}
					\toprule 
					\multicolumn{3}{c}{} & \multicolumn{2}{c}{} &  & \multicolumn{3}{l}{Model evaluations}\tabularnewline
					\multicolumn{3}{c}{Methods} & $\mathrm{CVaR}_{\beta}\,\text{estimate}{}^{(\mathrm{a})}$  & MRD$^{(\mathrm{a})}$ (\%)  & {N-RMSD$^{(\mathrm{a})}$ (\%)} & HF$^{(\mathrm{b})}$ & LF~\RNum{1}$^{(\rm{c})}$ & LF~\RNum{2}$^{(\rm{d})}$ \tabularnewline
					\midrule
					\multicolumn{9}{c}{{Correlation coefficient between $X_1$ and $X_2$ = $0.9$}}\tabularnewline
					\multicolumn{5}{l}{DD-GPCE-Kriging-based MCS} & \tabularnewline
					& \multicolumn{2}{c}{$S=1,\,m=1$} & $18.7363$ & $0.2046$ & $0.2667
					$ & $300$$^{(\mathrm{e})}$ & - & -\tabularnewline
					& \multicolumn{2}{c}{$S=1,\,m=2$} & $18.7366$ & $0.2042$ & $0.2640
					$ & $300$$^{(\mathrm{e})}$ & - & -\tabularnewline
					& \multicolumn{2}{c}{$S=1,\,m=3$} & $18.7366$ & $0.2043$ & $0.2641
					$ & $300$$^{(\mathrm{e})}$ & - & -\tabularnewline
					\multicolumn{5}{l}{DD-GPCE-based MCS} & \tabularnewline
					& \multicolumn{2}{c}{$S=1,\,m=1$} & $5.4175$ & $71.1383$ & $71.8302
					$ & $300$$^{(\mathrm{e})}$ & - & -\tabularnewline
					& \multicolumn{2}{c}{$S=1,\,m=2$} & $10.1197$ & $46.0874$ & $46.1142
					$ & $300$$^{(\mathrm{e})}$ & - & -\tabularnewline
					& \multicolumn{2}{c}{$S=1,\,m=3$} & $10.1499$ & $45.9263$ & $45.9534
					$ & $300$$^{(\mathrm{e})}$ & - & -\tabularnewline
					\multicolumn{5}{l}{PCE-Kriging-based MCS} &  &  &  & \tabularnewline
					& \multicolumn{2}{l}{$m=1$} & $18.0600$ & $3.9322$ & $4.3663
					$ & $300$$^{(\mathrm{e})}$ & - & -\tabularnewline
					& \multicolumn{2}{l}{$m=2$} & $18.0648$ & $3.9104$ & $4.3472
					$ & $300$$^{(\mathrm{e})}$ & - & -\tabularnewline
					& \multicolumn{2}{l}{$m=3$} & $18.0676$ & $3.8971$ & $4.3317
					$ & $300$$^{(\mathrm{e})}$ & - & -\tabularnewline
					\multicolumn{5}{l}{DD-GPCE-Kriging-based MFIS } & \tabularnewline
					& \multicolumn{2}{l}{HF} & $18.8607$ & $0.8361$ & $1.0674
					$ & $300$$^{(\mathrm{f})}$ & - & $10,000$$^{(\mathrm{g})}$\tabularnewline
					& \multicolumn{2}{l}{LF} &  &  &  &  &  & \tabularnewline
					&  & LF model \#1~\eqref{sur1} & $18.8217
					$ & $0.8119
					$ & $1.0545
					$& $150$$^{(\mathrm{h})}$ & $150$$^{(\mathrm{i})}$ & $10,000$$^{(\mathrm{g})}$\tabularnewline
					&  & LF model \#2~\eqref{sur2} & $18.8624
					$ & $0.8538
					$ & $1.0720
					$ & $150$$^{(\mathrm{h})}$ & $150$$^{(\mathrm{i})}$ & $10,000$$^{(\mathrm{g})}$\tabularnewline
					&  & LF model \#3~\eqref{sur3} & $16.0848
					$ & $14.3078
					$ & $14.5028
					$ & $150$$^{(\mathrm{h})}$ & $150$$^{(\mathrm{i})}$ & $10,000$$^{(\mathrm{g})}$\tabularnewline
					&  & LF model \#4~\eqref{sur4} & $18.2318
					$ & $2.8697
					$ & $3.0861
					$  & $150$$^{(\mathrm{h})}$ & $150$$^{(\mathrm{i})}$ & $10,000$$^{(\mathrm{g})}$\tabularnewline
					\multicolumn{3}{l}{Standard MCS (Benchmark)} & $18.7705$ & - &  & $10,000$ & - & -\tabularnewline
					\hdashline
					\multicolumn{9}{c}{ Correlation coefficient between $X_1$ and $X_2$ = $0$}\tabularnewline
					\multicolumn{5}{l}{DD-GPCE-Kriging-based MCS ($S=1,\,m=3$)} &  &  &  & \tabularnewline
					& \multicolumn{2}{l}{} & $18.3036$ & $0.4438$ & $ 0.6065$ & $300$$^{(\mathrm{e})}$ & - & -\tabularnewline
					\multicolumn{5}{l}{DD-GPCE-based MCS ($S=1,\,m=3$)} &  &  &  & \tabularnewline
					& \multicolumn{2}{l}{} & $10.0410$ & $45.3839$ & $45.4445$ & $300$$^{(\mathrm{e})}$ & - & -\tabularnewline
					\multicolumn{5}{l}{PCE-Kriging-based MCS ($m=3$)} &  &  &  & \tabularnewline
					& \multicolumn{2}{l}{} & $18.3007$ & $0.4623$ & $  0.6245$ & $300$$^{(\mathrm{e})}$ & - & -\tabularnewline
					\multicolumn{3}{l}{Standard MCS (Benchmark)} & $18.3848$ & - &  & $10,000$ & - & -\tabularnewline
					\bottomrule
			\end{tabular}}
			\par\end{center}
	\end{spacing}
	\vspace{-0.1in}
	\begin{tablenotes}
		\scriptsize\smallskip 
		\item{a.} The estimates are averaged over $K=50$ trials.
		\item{b.} The high-fidelity output is obtained by the Rastrigin function~\eqref{ras_f}.   
		\item{c.} The low-fidelity output is obtained by either of \eqref{sur1}--\eqref{sur4}. 
		\item{d.} The low-fidelity output is obtained by DD-GPCE-Kriging. 		
		\item{e.} The high-fidelity output samples are used to calculate DD-GPCE-Kriging, DD-GPCE, or PCE-Kriging. 
		\item{f.} Two different 50\% high-fidelity output samples are used to estimate the CI-based $\epsilon$-risk region and $\rm{CVaR}_{\beta}$, respectively.		
		\item{g.} The low-fidelity output samples are used to estimate the CI-based $\epsilon$ risk region.
		\item{h.} The high-fidelity output samples are used to estimate the $\rm{CVaR}_{\beta}$.
		\item{i.} The low-fidelity output samples are used to estimate DD-GPCE-Kriging.
	\end{tablenotes}
	\label{table1}
\end{table*}

Table~\ref{table1} presents the sampling-based CVaR estimates for the Rastrigin function \eqref{ras_f} when $\beta=0.99$. We show MCS estimates that use the DD-GPCE-Kriging, DD-GPCE, and PCE-Kriging as well as the high-fidelity model for comparison. We also show the DD-GPCE-Kriging-based MFIS estimator. 
The DD-GPCE-Kriging, DD-GPCE, and PCE-Kriging surrogates (whether using first, second or third-order approximation) are calculated using the same number (300) of high-fidelity output samples. We consider two separate cases where the two input variables have correlation coefficients 0.9 and 0, respectively. We provide benchmark estimates obtained by standard MCS with $10,000$ high-fidelity output samples.

When the inputs $X_1$ and $X_2$ are highly correlated with correlation coefficient 0.9, the univariate ($S=1$) DD-GPCE-Kriging methods of first ($m=1$) through third ($m=3$) order approximations provide MCS-based CVaR estimates that are close to the standard MCS. Their MRD over $K=50$ trials are around $0.20$\%. In contrast, the MCS-based CVaR estimates by univariate ($S=1$) DD-GPCE or PCE-Kriging methods of the same order ($m=1-3$) approximations are less accurate compared with the standard MCS solution. The accuracy of CVaR estimates by the DD-GPCE or PCE-Kriging method improves from $71.14\%$ to $45.92\%$ or from $3.93\%$ to $3.90\%$, respectively, in MRD as the order of approximation increases from $m=1$ to $m=3$. For the highly nonlinear Rastrigin function, the MCS-based CVaR estimate that uses the DD-GPCE-Kriging surrogate is therefore significantly more accurate than the estimate by MCS that samples from the DD-GPCE and PCE-Kriging surrogates. 

When the inputs $X_1$ and $X_2$ have zero correlation, the accuracy of the proposed DD-GPCE-Kriging-based MCS method is still better than the DD-GPCE method. However, the PCE-Kriging method performs similar to DD-GPCE-Kriging, showing only a $0.02\%$ difference between their MRD values. The N-RMSD values of DD-GPCE-Kriging and PCE-Kriging methods are both around $0.6\%$, indicating almost the same standard deviation of both estimators relative to the benchmark. This is because both DD-GPCE-Kriging and PCE-Kriging use the same orthonormal basis functions that are measure-consistent with independent random inputs.
However, for dependent inputs, only DD-GPCE-Kriging can use such basis functions without any transformation from dependent to independent random inputs. Therefore, the DD-GPCE-Kriging outperforms PCE-Kriging in prediction accuracy when the random input variables are dependent. We also conclude that the additional variance term introduced in DD-GPCE-Kriging can better capture the local oscillations in the Rastrigin function compared to the more global approximation that is DD-GPCE.  

We next present results for the proposed MFIS-based method. Recall, that option ``HF'' uses the high-fidelity model to learn the DD-GPCE-Kriging surrogate, and option ``LF'' uses the low-fidelity model to learn the DD-GPCE-Kriging surrogate. In both cases, those surrogates are only used to learn the biasing density for importance sampling. The univariate ($S=1$), third-order ($m=3$) DD-GPCE-Kriging was chosen to determine the CI-based $\epsilon$-risk region for all MFIS cases tabulated in the second through sixth rows from the bottom of Table~\ref{table1}. The proposed MFIS-based CVaR estimate is 18.42 with HF, using 300 high-fidelity output samples. That estimate is close to the benchmark estimate of 18.38 by standard MCS. However, its MRD value of 0.80\% is slightly higher than the 0.36--0.37\% computed for the DD-GPCE-Kriging methods that also use 300 high-fidelity output samples. To speed up the process of learning the biasing density for MFIS, we consider the LF option. Specifically, in this example, we are interested in the accuracy of the MFIS-based CVaR estimates by LF when each of the four distinct surrogate models~\eqref{sur1}--\eqref{sur4} is used. The LF option with low-fidelity models \#1 and \#2 ($\rho_{H,L_1}$ and $\rho_{H,L_2}=1$) yield a similar accuracy to the HF option. Contrarily, the LF option with low-fidelity models \#3 and \#4, for which $\rho_{H,L_3}$ and $\rho_{H,L_4}=0.72$, provides less accurate CVaR estimates. The results indicate that for the proposed MFIS-based method, it is acceptable for the low-fidelity models ($y_{L_i}$,~$i$=1--4) to have a bias (which low-fidelity models~\#1 and \#2 certainty do) as long as they have a strong correlation with the high-fidelity output $y_H$.  

\subsection{Example 2: Modified cross-in-tray function} \label{sec:5.4}
This second example illustrates the DD-GPCE-Kriging-based MCS and MFIS methods for a nonsmooth analytical mathematical function. 

\subsubsection{Problem definition}

Define a bivariate Gaussian random vector $\bX=(X_1,X_2)^{\intercal}$ with mean vector $(0,0)^{\intercal}$ and standard deviations $\sqrt{\mathrm{var}[X_i]}=2,~i=1,2$. The correlation coefficient between $X_1$ and $X_2$ is $0.9$. For the input vector $\bX$, consider the modified cross-in-tray function:

\begin{align} \label{mod_f} 
	y(X_1,X_2)=-0.001\left(\bigg|\sin(X_1) \sin(X_2) \exp{\left(\bigg|100-\sqrt{\dfrac{X_1^2+X_2^2}{\pi}}\bigg|\right)}\bigg|+1\right)^{0.1}.
\end{align} 

Figure~\ref{fig2}b shows the modified cross-in-tray function on the domain $[-4,4]\times[-4,4]$. The function has two ridges along the $x_1=0$ and $x_2=0$ axes, respectively. It sharply drops on each side of the ridges, forming sharp corners and exhibiting oscillatory behavior. The sharp cross is not differentiable. The function is notoriously difficult to be approximate with polynomial bases due to its nonsmoothness.  

\subsubsection{Results}
\begin{table*}
	\caption{$\mathrm{CVaR}_{\beta}$ estimates ($\beta=0.99$) of the modified cross-in-tray function~\eqref{mod_f} via the DD-GPCE-Kriging, DD-GPCE, and PCE-Kriging-based MCS methods, the DD-GPCE-Kriging-based MFIS method, and standard MCS.}
	\vspace{-0.1in}
	\begin{spacing}{1.2}
		\begin{center}
			\small
			\resizebox{\textwidth}{!}{
				\begin{tabular}{clcccccc}
					\toprule 
					\multicolumn{3}{c}{} & \multicolumn{2}{c}{} &  & \multicolumn{2}{c}{Model evaluations}\tabularnewline
					\multicolumn{3}{c}{Methods} & $\mathrm{CVaR}_{\beta}\,\text{estimate}{}^{(\mathrm{a})}$  & MRD$^{(\mathrm{a})}$ (\%)  & {N-RMSD$^{(\mathrm{a})}$ (\%)} & HF$^{(\mathrm{b})}$ & LF$^{(\mathrm{c})}$\tabularnewline
					\midrule
					\multicolumn{5}{l}{DD-GPCE-Kriging-based MCS} & \tabularnewline
					& \multicolumn{2}{c}{$S=1,\,m=3$} & $-10.4649$ & $6.1486$ & $6.8785$ & $400$$^{(\mathrm{d})}$ & -\tabularnewline
					& \multicolumn{2}{c}{$S=1,\,m=4$} & $-10.4040$ & $5.6964$ & $6.4619$ & $400$$^{(\mathrm{d})}$ & -\tabularnewline
					& \multicolumn{2}{c}{$S=1,\,m=5$} & $-10.3035$ & $4.9523$ & $5.7967$ & $400$$^{(\mathrm{d})}$ & -\tabularnewline
					\multicolumn{5}{l}{DD-GPCE-based MCS} & \tabularnewline
					& \multicolumn{2}{c}{$S=1,\,m=3$} & $-13.4384$ & $35.9874$ & $36.7051$ & $400$$^{(\mathrm{d})}$ & -\tabularnewline
					& \multicolumn{2}{c}{$S=1,\,m=4$} & $-14.2902$ & $44.6070$ & $45.5094$ & $400$$^{(\mathrm{d})}$ & -\tabularnewline
					& \multicolumn{2}{c}{$S=1,\,m=5$} & $-13.9834$ & $41.5772$ & $43.5414$ & $400$$^{(\mathrm{d})}$ & -\tabularnewline
					\multicolumn{5}{l}{PCE-Kriging-based MCS} &  &  & \tabularnewline
					& \multicolumn{2}{c}{$m=3$} & $-11.1909$ & $13.2435$ & $14.7770$ & $400$$^{(\mathrm{d})}$ & -\tabularnewline
					& \multicolumn{2}{c}{$m=4$} & $-10.8560$ & $12.6933$ & $14.8892$ & $400$$^{(\mathrm{d})}$ & -\tabularnewline
					& \multicolumn{2}{c}{$m=5$} & $-10.1499$ & $12.4506$ & $19.6315$ & $400$$^{(\mathrm{d})}$ & -\tabularnewline
					\multicolumn{5}{l}{DD-GPCE-Kriging-based MFIS (HF)} & \tabularnewline
					& \multicolumn{2}{c}{$m=3$} & $-10.2050$ & $3.3645$ & $3.9706$ & $400$$^{(\mathrm{e})}$ & $10,000$$^{(\mathrm{f})}$\tabularnewline
					& \multicolumn{2}{c}{$m=4$} & $-10.1884$ & $3.1556$ & $3.6153$ & $400$$^{(\mathrm{e})}$ & $10,000$$^{(\mathrm{f})}$\tabularnewline
					& \multicolumn{2}{c}{$m=5$} & $-10.2141$ & $3.3589$ & $3.7319$ & $400$$^{(\mathrm{e})}$ & $10,000$$^{(\mathrm{f})}$\tabularnewline
					\multicolumn{3}{l}{Standard MCS (Benchmark)} & $-9.8821$ & - &  & $10,000$ & -\tabularnewline
					\bottomrule
			\end{tabular}}
			\par\end{center}
	\end{spacing}
	\vspace{-0.1in}
	\begin{tablenotes}
		\scriptsize\smallskip 
		\item{a.} The estimates are averaged over $K=50$ trials.
		\item{b.} The high-fidelity output is obtained by the modified cross-in-tray function~\eqref{mod_f}.    
		\item{c.} The low-fidelity output is obtained by the DD-GPCE-Kriging. 	
		\item{d.} The high-fidelity output samples are used to calculate DD-GPCE-Kriging, DD-GPCE, or PCE-Kriging. 
		\item{e.} Two different 50\% high-fidelity output samples are used to estimate the CI-based $\epsilon$-risk region and $\rm{CVaR}_{\beta}$, respectively.	
		\item{f.} The low-fidelity output samples are used to estimate the CI-based $\epsilon$-risk region. 
	\end{tablenotes}
	\label{table1:2}
\end{table*}

Table~\ref{table1:2} presents the sampling-based CVaR estimates of the modified cross-in-tray function \eqref{mod_f} when $\beta=0.99$. We show MCS estimates that use the DD-GPCE-Kriging, DD-GPCE, and PCE-Kriging; we provide a benchmark estimate obtained by standard MCS of the high-fidelity model with $10,000$ samples for comparison. We also show the DD-GPCE-Kriging-based MFIS estimator. 

We evaluate the accuracy of the DD-GPCE-Kriging-based MCS and MFIS methods when the quantity of interest includes nonsmooth features, as shown in Figure~\ref{fig2}b. While the DD-GPCE-Kriging-based MCS improves the accuracy of CVaR estimation as the order increases from $m=3$ to $m=5$, the MFIS-based CVaR estimates with the HF option are the most accurate with 3.16--3.36\% MRD and are relatively insensitive to the order ($m$) increase of the DD-GPCE-Kriging. This is because for the CVaR estimate, MFIS uses high-fidelity output samples drawn from the biasing density. Here, we only use a small number of 200 high-fidelity output samples for CVaR estimation to demonstrate its efficiency. While the DD-GPCE-Kriging surrogate is somewhat limited in predicting the nonsmooth output, the surrogate is only used to learn the biasing density. The standard deviation of the MFIS estimators relative to the benchmark is also the smallest compared to the MCS-based methods, in the range of $3.6$--$3.97\%$ N-RMSD. In sum, for the nonsmooth modified cross-in-tray function, the proposed MFIS-based method is more accurate in providing CVaR estimates than DD-GPCE-Kriging-based MCS.  

\subsubsection{Convergence properties of the DD-GPCE-Kriging method}
		\begin{figure*}
			\begin{center}
				\includegraphics[angle=0,scale=0.9,clip]{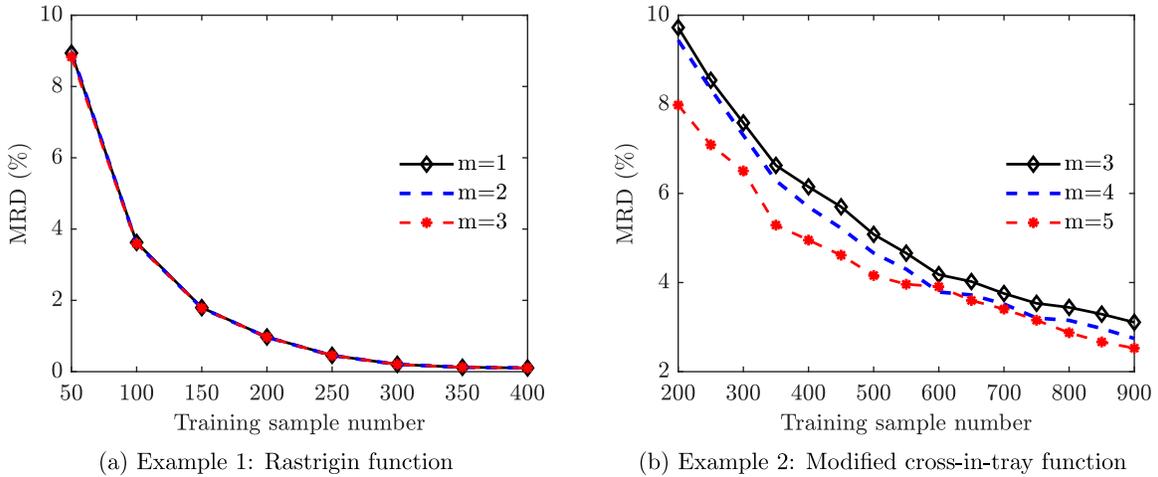}
			\end{center}
			\caption{{MRD (\%) versus training data size when estimating CVaR by Monte Carlos sampling of the DD-GPCE-Kriging surrogates. Various degrees ($m$) of the polynomial approximation are illustrated.}}
			\label{fig3}
		\end{figure*} 
		
		We examine the convergence of the DD-GPCE-Kriging-based Monte Carlo estimator with respect to the sample size that is used for training the DD-GPCE-Kriging surrogate. We illustrate those on Examples 1 and 2. 
		Figures~\ref{fig3}(a) and (b) show MRD versus training sample size when the CVaR estimates are obtained by sampling directly from the DD-GPCE-Kriging models. The MRD is computed relative to the MCS benchmark in Example 1 and 2. 
		In Example~1, the MRD decreases from 8.8\% to 0.1\% as the training sample number increases from a small number of 50 to a modest 400. The behavior is almost identical for the three DD-GPCE-Kriging surrogates with different order approximations ($m=1-3$). The figure illustrates that the DD-GPCE-Kriging method is---as expected in such high-dimensional problems---less accurate when a small number of samples are used to train the surrogate. 
		When considering the results for Example~2 in Figure~\ref{fig3}(b), we observe a similar monotone decrease in MRD as the sample number increases. This time, however, there is also a noticeable improvement from the $m=3$ to $m=4$ and $m=5$ DD-GPCE-Kriging approximations. 
		In both cases of heavily nonlinear and/or nonsmooth functions, and for all $m$, the DD-GPCE-Kriging method becomes more accurate as the training sample number increases.

\subsection{Example 3: A 2D glass/vinylester composite plate} \label{sec:5.5}

This example studies a bolted or riveted composite laminate, which is widely used in aircraft structures. We are interested in measuring the risk of failure that occurs by stress concentration near rivet holes. The quantity of interest is obtained by solving a nonlinear quasi-static two-dimensional FEA. The model has a high number ($N=28$) of random inputs, some of which are modeled as dependent random variables. This provides a challenging case for risk assessment. 
\subsubsection{Problem description}

\begin{figure*}
	\begin{center}
		\includegraphics[angle=0,scale=0.9,clip]{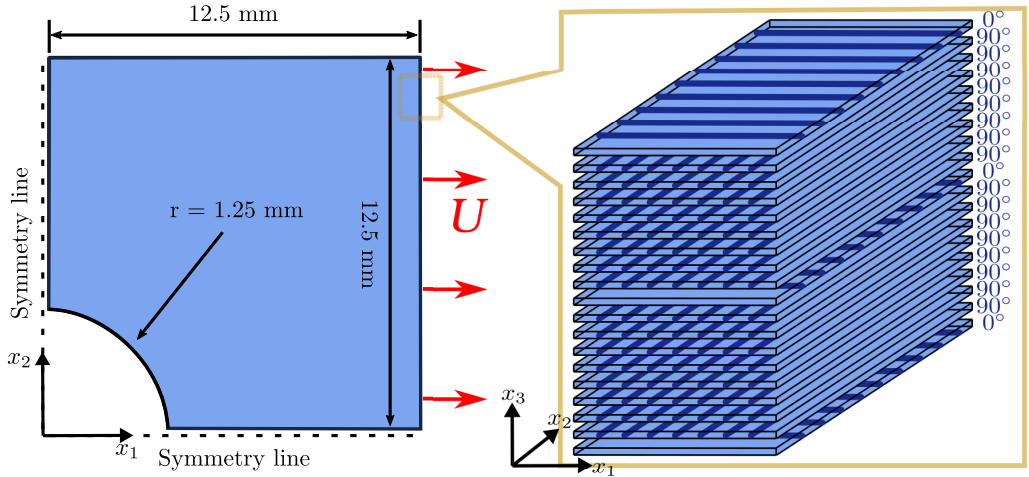}
	\end{center}
	\vspace{-0.2in}
	\caption{Geometry, loading, and boundary conditions of a glass/vinylester composite plate: We consider one quarter of the full plate domain, which represents the original square plate with domain $\mathcal{D}=(25~\mathrm{mm}\times25~\mathrm{mm})$ including a circular hole of radius $r=1.25~\mathrm{mm}$ in the center due to the symmetry conditions in the $x_1$ and $x_2$ directions of the center lines.}
	\label{fig4}
\end{figure*}		
Figure~\ref{fig4}a presents a quarter model of a two-dimensional square plate (spatial domain $\mathcal{D}=25~\mathrm{mm}\times25~\mathrm{mm}$) with a circular hole of radius $r=1.25~\mathrm{mm}$ in the center. Due to symmetry conditions with respect to $x_1$ and $x_2$, we can consider the quarter model in the domain $\bar{\cD}=(12.5~\mathrm{mm}\times12.5~\mathrm{mm})$ to save computational cost. 
Figure~\ref{fig4}b shows the arrangement of the quarter model made of the glass/vinylester laminate (Fiberite/HyE 9082Af) with 19 stacking sequences $[0/90_8/0/90_8/0]$, where `0' indicates a ply having fibers in $x_1$ direction and $90_8$ denotes that eight plies have fibers in $x_2$ direction.
The plate is subjected to a uniaxial tensile displacement loading $U$ (mm/s) that acts on the entire right side edge.
This example can also be found in~\cite{MOURE2014224,lee2022bi}.  
\begin{table*}
	\caption{Properties of the random inputs in Fiberite/HyE 9082Af}
	\begin{spacing}{1.0}
		\begin{centering}
			\footnotesize				
			\begin{tabular}{ccccccc}
				\toprule 
				Random  & \multirow{2}{*}{Property} & \multirow{2}{*}{Mean} & {Coefficient of variation} & Lower  & Upper & Probability \tabularnewline
				variable &  &  & {($\rm{COV}_i^{(\rm{a})}$, \%)} & boundary & boundary & distribution\tabularnewline
				\midrule 
				$X_{1}$ & $E_{1}$ (MPa) & $44,700$ & 11.55  & 35760 & 53640 & Uniform\tabularnewline
				$X_{2}$ & $E_{2}$ (MPa) & $12,700$ & 11.55  & 10,160 & 15,240 & Uniform\tabularnewline
				$X_{3}$ & $v_{12}$ & $0.297$ & 11.55  & 0.238 & 0.356 & Uniform\tabularnewline
				$X_{4}$ & $G_{12}$ (MPa) & $5,800$ & 11.55  & 4,640 & 6,960 & Uniform\tabularnewline
				$X_{5}$ & $S_{t1}$ (MPa) & $1,020$ & 11.55  & 816 & 1,224 & Uniform\tabularnewline
				$X_{6}$ & $S_{t2}$ (MPa) & $40$ & 11.55  & 32 & 48 & Uniform\tabularnewline
				$X_{7}$ & $S_{c1}$ (MPa) & $620$ & 11.55  & 496 & 744 & Uniform\tabularnewline
				$X_{8}$ & $S_{c2}$ (MPa) & $140$ & 11.55  & 126 & 168 & Uniform\tabularnewline
				$X_{9}$ & $S_{s12}$$^{(\mathrm{b})}$ (MPa) & $60$ & 11.55  & 48 & 72 & Uniform\tabularnewline
				\multirow{2}{*}{$X_{10}$--$X_{28}$$^{(\mathrm{c})}$} & Plies 1--19 & \multirow{2}{*}{$0.144$} & \multirow{2}{*}{$6$} & \multirow{2}{*}{0} & \multirow{2}{*}{$\infty$} & Multivariate\tabularnewline
				& thicknesses (mm) &  &  &  &  & Lognormal\tabularnewline
				\bottomrule
			\end{tabular}
			\par\end{centering}
	\end{spacing}\tabularnewline
	\begin{tablenotes}
		\scriptsize\smallskip
		\item{a.} $\rm{COV}_i=100\times\sqrt{\rm{var}[X_i]}/\mathbb{E}{[X_i]}$, $i=1,\ldots,28$.
		\item{b.} $S_{s12}=S_{s23}.$
		\item{c.} Correlation coefficients among $X_{10}$--$X_{28}$ are 0.5.
	\end{tablenotes}			
	\label{table2}
\end{table*}
The local spatial coordinates are $x_1,~x_2,~\text{and}~x_3$, and we define $E_i$, $\nu_{ij}$, and $G_{ij}$ to be Young's modulus, Poisson's ratio, and shear modulus of the plies in the corresponding local spatial coordinates for $i,j=1,2,3$. Let $S_{ti}$, $S_{ci}$, and $S_{sij}$ be the tensile, compressive, and shear strengths of the plies. The measured material properties and ply thicknesses of the laminate vary so we model them as $N=28$ random variables, as presented in Table~\ref{table2}. The $19$ random variables that model the ply thickness are modeled as correlated via a multivariate lognormal distribution with a correlation coefficient of $0.5$. The remaining $9$ random variables are modeled as independent and distributed uniformly.

\subsubsection{Quantity of interest}

\begin{figure*}
	\begin{center}
		\includegraphics[angle=0,scale=0.8,clip]{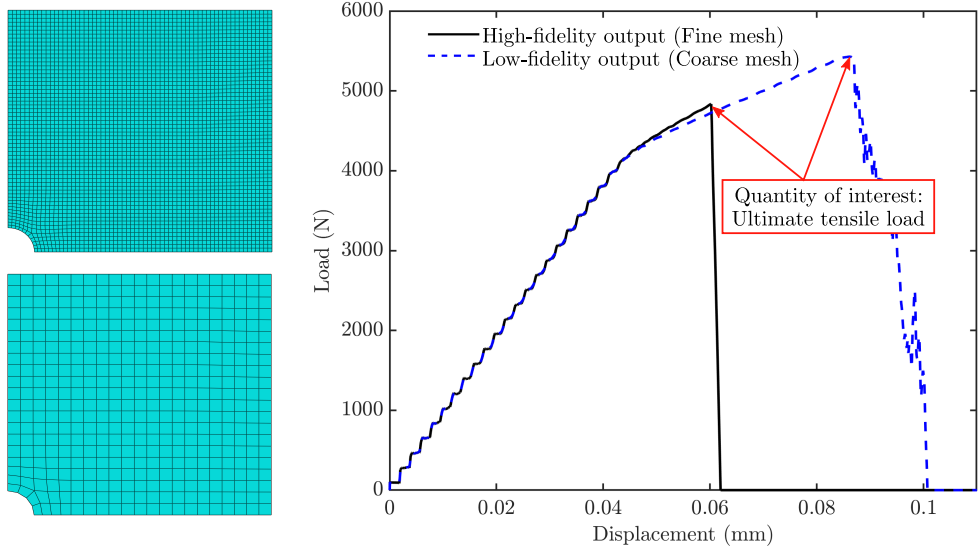}
	\end{center}
	\caption{Comparison of meshes and the load-displacement curve between high-fidelity and low-fidelity models: Figure~(a) presents the fine mesh with 3,887 elements (Top) and the coarse mesh with 441 elements (Bottom). The fine and coarse mesh models generate high-fidelity and low-fidelity output data. In plot~(b), high-fidelity and low-fidelity outputs in the load-displacement curve have different maximum values (ultimate tensile load).}
	\label{fig5}
\end{figure*}	
The high-fidelity model is a finite element model with a fine mesh, and the low-fidelity model is a finite element model with a coarse mesh, see Figure~\ref{fig5}a. For both models, the element types were chosen as S4R, a 4-node, quadrilateral, stress/displacement shell element with reduced integration, built-in ABAQUS/Explicit, version 6.14-2. The total number of degrees of freedom for the high-fidelity and low-fidelity models are $24,084$ and $2,910$, respectively. 

Figure~\ref{fig5}b shows the load versus displacement curves of both the high-fidelity and low-fidelity models. Given the increasing applied displacement $U$ in the range of $0$--$0.04$~$\mathrm{mm}$ at the right edges of both models, the resultant tensile loads increase almost linearly in each case. After a displacement of about $0.04$mm, the rate of change in both load curves decreases. The load values peak right before fracture. We then see a failure/drop right after the maximum tensile load is reached. We consider the peak load values, referred to as an ultimate tensile load, as the output of interest.  The damage emanates from the hole and spreads outwards to the plate until it is completely fractured. To describe the damage, we use the Hashin damage model~\cite{hashin1980} built-in ABAQUS/Explicit, version 6.14-2, see \cite{lee2022bi} for a detailed explanation. 

The correlation coefficient between the high-fidelity model $\rho_{H,L}$ and the low-fidelity model $y_L$ is $\rho_{H,L}=0.80$, see \eqref{pcor} for its computation. This indicates that the low-fidelity and high-fidelity model are strongly correlated. Thus, we use the low-fidelity model to compute the DD-GPCE-Kriging surrogate, which is then used to learn the biasing density for MFIS.  

\subsubsection{Results}
\begin{table*}
	\caption{$\mathrm{CVaR}_{\beta}$ estimates ($\beta=0.99$) of the ultimate tensile load of the glass/vinylester composite plate via the DD-GPCE-Kriging, DD-GPCE, and PCE-Kriging-based MCS methods, DD-GPCE-Kriging-based MFIS, and standard MCS.}
	\vspace{-0.1in}
	\begin{spacing}{1.2}
		\begin{center}
			\small
			\resizebox{\textwidth}{!}{
				\begin{tabular}{clcccccccc}
					\toprule 
					\multicolumn{3}{c}{} & \multicolumn{3}{l}{Ultimate tensile load (N)} & \multicolumn{3}{l}{Model evaluations} & CPU time\tabularnewline
					\multicolumn{3}{c}{Methods} & $\mathrm{CVaR}_{\beta}$ estimate$^{(\mathrm{a})}$  & MRD (\%)$^{(\mathrm{a})}$  & {N-RMSD}$^{(\mathrm{a})} (\%)$  & HF$^{(\mathrm{b})}$ & LF~\RNum{1}$^{(\mathrm{c})}$ & LF~\RNum{2}$^{(\mathrm{d})}$ & (hours)$^{(\mathrm{e})}$\tabularnewline
					\midrule
					\multicolumn{5}{l}{DD-GPCE-Kriging-based MCS} &  & \tabularnewline
					& \multicolumn{2}{c}{$S=1,\,m=1$} & $5558.0651$ & $0.4383$ & $0.4985$ & $250$$^{(\mathrm{f})}$ & - & - & $21.5$\tabularnewline
					& \multicolumn{2}{c}{$S=1,\,m=2$} & $5556.7570$ & $0.4861$ & $0.6189$ & $250$$^{(\mathrm{f})}$ & - & - & $21.5$\tabularnewline
					& \multicolumn{2}{c}{$S=1,\,m=3$} & $5556.7570$ & $0.4861$ & $0.5815$ & $250$$^{(\mathrm{f})}$ & - & - & $21.5$\tabularnewline
					\multicolumn{5}{l}{DD-GPCE-based MCS} &  & \tabularnewline
					& \multicolumn{2}{c}{$S=1,\,m=1$} & $5529.3800$ & $0.9381$ & $0.9570$ & $250$$^{(\mathrm{f})}$ & - & - & $21.5$\tabularnewline
					& \multicolumn{2}{c}{$S=1,\,m=2$} & $5533.2096$ & $0.8743$ & $0.9164$ & $250$$^{(\mathrm{f})}$ & - & - & $21.5$\tabularnewline
					& \multicolumn{2}{c}{$S=1,\,m=3$} & $5542.4489$ & $0.7334$ & $0.8220$ & $250$$^{(\mathrm{f})}$ & - & - & $21.5$\tabularnewline
					\multicolumn{5}{l}{DD-GPCE-Kriging-based MFIS } &  & \tabularnewline
					& \multicolumn{2}{l}{HF} & $5577.0842$ & $0.5441$ & $0.6784$ & $250$$^{(\mathrm{g})}$ & - & $10,000$$^{(\mathrm{h})}$ & $21.5$\tabularnewline
					& \multicolumn{2}{l}{LF} & $5577.1421$ & $1.1459$ & $1.6112$ & $50$$^{(\mathrm{i})}$ & $200$$^{(\mathrm{j})}$ & $10,000$$^{(\mathrm{h})}$ & $8.2$\tabularnewline
					\multicolumn{3}{l}{Standard MCS (Benchmark)} & $5582.1484$ & - &  & $10,000$ & - & - & $859.8$\tabularnewline
					\bottomrule
			\end{tabular}}
			\par\end{center}
	\end{spacing}
	\vspace{-0.1in}
	\begin{tablenotes}
		\scriptsize\smallskip 
		\item{a.} The estimates are averaged over $K=30$ trials.
		\item{b.} The high-fidelity output is obtained by the fine mesh model in Figure~\ref{fig5}a~(Top).   
		\item{c.} The low-fidelity output is obtained by the coarse mesh model in Figure~\ref{fig5}a~(Bottom).  
		\item{d.} The low-fidelity output is obtained by the DD-GPCE-Kriging. 	
		\item{e.} The total CPU time is (the number of FEA)$\times$(FEA runtime) in each trial.
		\item{f.} The high-fidelity output samples are used to calculate the DD-GPCE-Kriging or DD-GPCE.  
		\item{g.} The 250 samples are split as follows: 200 samples are used to compute DD-GPCE-Kriging and 50 samples are used to estimate~$\rm{CVaR}_{\beta}$.   
		\item{h.} The low-fidelity output samples are used to estimate the CI-based $\epsilon$-risk region.
		\item{i.} The high-fidelity output samples are used to estimate the $\rm{CVaR}_{\beta}$.
		\item{j.} The low-fidelity output samples are used to estimate the DD-GPCE-Kriging. 		
	\end{tablenotes}
	\label{table3}
\end{table*}

Table~\ref{table3} presents the CVaR estimates obtained from the DD-GPCE-Kriging and DD-GPCE-based MCS methods and DD-GPCE-Kriging-based MFIS methods with two options (Option ``LF" trains the DD-GPCE-Kriging with the low-fidelity model and Option ``HF" with the high-fidelity model). For comparison, we also provide a benchmark estimate obtained by the standard MCS with 10,000 samples in the last row of Table~\ref{table3}. The CVaR estimates are accurate with $0.43\%$ through $1.15\%$ in MRD over $K=30$ trials. The CVaR estimates by the univariate ($S=1$) DD-GPCE-Kriging methods of the first through third-order ($m=1-3$) approximations are more precise (about two times lower MRD) than those by DD-GPCE when the same number ($300$) of high-fidelity output samples are used. 
The MFIS with two options (HF and LF) exhibits good accuracy (0.54\% and 1.15\%, respectively in MRD) of the CVaR estimation. The MFIS with LF has an N-RMSD of $1.6\%$, which is relatively large compared to those ($0.5\%$--$0.6\%$) of MCS-based methods or the MFIS method with HF. The LF option requires a CPU time of 8.2 hours, which is an almost 62\% reduction compared to those (21.5 hours) required by the HF option. 
Such computational gains are made feasible by using the cheaper low-fidelity surrogate for constructing the DD-GPCE-Kriging model.

\subsection{Example 4: A 3D composite T-joint} \label{sec:5.6}
Our final example considers a three-dimensional carbon fiber-reinforced silicon carbide (C/SiC) T-joint that requires a nonlinear quasi-static FEA to obtain the quantity of interest. This model has $N=20$ random inputs, some of which are dependent. The complexity of the 3D FEA and the high-dimensional inputs make this a challenging case for risk analysis.   
\begin{figure*}
	\begin{center}
		\includegraphics[angle=0,scale=0.8,clip]{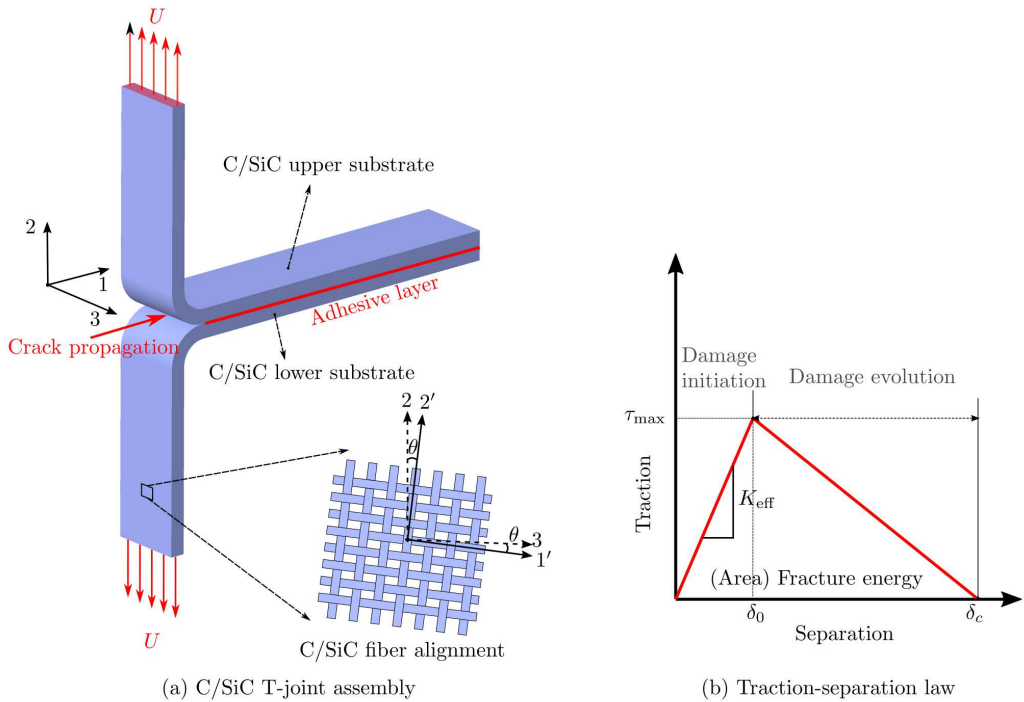}
	\end{center}
	\vspace{-0.2in}
	\caption{Geometry, loading, and boundary conditions of the 3D C/SiC composite T-joint assembly: (a) The two substrates consisting of interlaced C/SiC fibers are jointed by an adhesive layer; (b) The traction-separation law determines the damage of the adhesive layer under the peel loading $U$, leading to the deposition of the C/SiC T-joint. In Figure~\ref{fig6}b, $\tau_{\max}$, $K_{\rm{eff}}$, and the fracture energies are all considered to be uncertain, see Table~\ref{table4}.}
	\label{fig6}
\end{figure*}
\subsubsection{Problem description} \label{sec:5.4:1}

Figure~\ref{fig6}a shows a three-dimensional T-joint of interlaced C/SiC composite materials. It also illustrates the interlaced fiber alignments in the upper and lower substrates.  Since the composite consists of high-strength carbon fibers and a high modulus, oxidation-resistant matrix, it exhibits high specific strength, low coefficient of thermal expansion, moderate thermal conductivity, and is lightweight. This makes this composite material appealing for super and hyper-sonic transport propulsion systems (e.g., space- or aircraft). The T-joint assembly allows the assembled components to avoid stress concentration occurring in the region near rivet holes and to achieve a lightweight assembled structure, see \cite{yan2018experimental}. 

The T-joint structure is subjected to a uni-axial tensile displacement loading $U$ that acts on both the top and bottom sides of the edges (known as a peel test). The peel test evaluates the performance of interfacial adhesion in the T-joint assembly. We describe interfacial adhesion through cohesive surfaces with the traction-separation law illustrated in Figure~\ref{fig6}b. In the graph of Figure~\ref{fig6}b, the traction and separation represent the nominal stress and the crack opening displacement, respectively, which occur at the crack front line of the adhesive layer. The traction-separation model assumes initially linear elastic behavior of the adhesive layer. After the nominal traction stress approaches $\tau_{\max}$, the layer damage initiates at the separation $\delta_0$. Then, the damage propagates until complete separation occurs at $\delta_c$. Given $\tau_{\max}$ and $K_{\rm{eff}}$, we determine the separation or crack opening displacement by the fracture energy (i.e., the integral of the traction over the separation in Figure~\ref{fig6}b).       

The measured material properties of C/SiC and adhesive parameters vary. Thus, we model them as $N=20$ random variables $X_i$, $i=1,\ldots,20$, following either multivariate lognormal or uniform distributions, see Table~\ref{table4}. For the material properties of the C/SiC composite, $E_{ij}$, $\nu_{ij}$, and $G_{ij}$ indicate elastic modulus, poison ratio, and shear modulus, respectively, in the local $j$-direction on the surface normal to the local $i$-direction, $i,j=1,2,3$. For the cohesive parameters associated with the T-joint deposition (detachment between two substrates), $K_{ii}$ and $\tau_{ii,\max}$ indicate the stiffness of the adhesive layer in the local $i$-direction, $i=1,2,3$, and the components of the maximum nominal traction $\tau_{\max}$, see Figure~\ref{fig6}b. The mean values of $X_i$, $i=1,\ldots,20$, were obtained from \cite{wang2021experimental}.    
\begin{table*}
	\caption{Properties of the random inputs for the C/SiC material and the adhesive layer.}
	\begin{spacing}{1.2}
		\begin{centering}
			\footnotesize	
			\begin{tabular}{ccccccc}
				\toprule 
				Random  & \multirow{2}{*}{Property} & \multirow{2}{*}{Mean} & {Coefficient of variation} & Lower  & Upper & Probability \tabularnewline
				variable &  &  & {($\rm{COV}_i^{\rm{(a)}}$, \%)}  & boundary & boundary & distribution\tabularnewline
				\midrule\tabularnewline
				\multicolumn{7}{c}{Material properties of C/SiC composite}\tabularnewline\tabularnewline
				$X_{1}$ & $E_{11}$(MPa) & $92,140$ & $5$ & 0 & $\infty$ & Lognormal$^{\rm{(b)}}$\tabularnewline
				$X_{2}$ & $E_{22}$(MPa) & $92,140$ & $5$ & 0 & $\infty$ & Lognormal$^{\rm{(b)}}$\tabularnewline
				$X_{3}$ & $E_{33}$(MPa) & $40,000$ & $5$ & 0 & $\infty$ & Lognormal$^{\rm{(b)}}$\tabularnewline
				$X_{4}$ & $v_{12}$ & $0.01$ & $5$ & 0 & $\infty$ & Lognormal$^{\rm{(b)}}$\tabularnewline
				$X_{5}$ & $v_{13}$ & $0.01$ & $5$ & 0 & $\infty$ & Lognormal$^{\rm{(b)}}$\tabularnewline
				$X_{6}$ & $v_{23}$ & $0.01$ & $5$ & 0 & $\infty$ & Lognormal$^{\rm{(b)}}$\tabularnewline
				$X_{7}$ & $G_{12}$(MPa)  & $23,180$ & $5$ & 0 & $\infty$ & Lognormal$^{\rm{(b)}}$\tabularnewline
				$X_{8}$ & $G_{13}$(MPa) & $23,180$ & $5$ & 0 & $\infty$ & Lognormal$^{\rm{(b)}}$\tabularnewline
				$X_{9}$ & $G_{23}$(MPa) & $23,180$ & $5$ & 0 & $\infty$ & Lognormal$^{\rm{(b)}}$\tabularnewline
				\tabularnewline\multicolumn{7}{c}{Cohesive parameters for the T-joint deposition}\tabularnewline\tabularnewline
				$X_{10}$ & $K_{11}$(MPa/mm) & $679,800$ & $10$ & $611,820$ & $747,780$ & Uniform\tabularnewline
				$X_{11}$ & $K_{22}$(MPa/mm) & $283,240$ & $10$ & $254,916$ & $311,564$ & Uniform\tabularnewline
				$X_{12}$ & $K_{33}$(MPa/mm) & $283,240$ & $10$ & $254,916$ & $311,564$ & Uniform\tabularnewline
				$X_{13}$ & $\tau_{11,\max}$(MPa) & $9.5$ & $10$ & $8.55$ & $10.45$ & Uniform\tabularnewline
				$X_{14}$ & $\tau_{22,\max}$(MPa) & $8.8$ & $10$ & $7.92$ & $9.68$ & Uniform\tabularnewline
				$X_{15}$ & $\tau_{33,\max}$(MPa) & $8.8$ & $10$ & $7.92$ & $9.68$ & Uniform\tabularnewline
				$X_{16}$ & Normal fracture energy(mJ) & $0.15$ & $10$ & $0.135$ & $0.165$ & Uniform\tabularnewline
				$X_{17}$ & $1$st shear fracture energy(mJ) & $0.15$ & $10$ & $0.135$ & $0.165$ & Uniform\tabularnewline
				$X_{18}$ & 2nd shear fracture energy(mJ) & $0.15$ & $10$ & $0.135$ & $0.165$ & Uniform\tabularnewline
				\multirow{2}{*}{$X_{19}$} & Fiber alignment in& \multirow{2}{*}{$0$} & \multirow{2}{*}{-$^{\rm{(c)}}$} & \multirow{2}{*}{$-10$} & \multirow{2}{*}{$10$} & \multirow{2}{*}{Uniform}\tabularnewline
				&upper substrate(degree) &  &  &  &  & \tabularnewline
				\multirow{2}{*}{$X_{20}$} & Fiber alignment in& \multirow{2}{*}{$0$} & \multirow{2}{*}{-$^{\rm{(c)}}$} & \multirow{2}{*}{$-10$} & \multirow{2}{*}{$10$} & \multirow{2}{*}{Uniform}\tabularnewline
				&lower substrate(degree) &  &  &  &  & \tabularnewline
				\bottomrule
			\end{tabular}
			\par\end{centering}
	\end{spacing}\tabularnewline
	\begin{tablenotes}
		\scriptsize\smallskip
		\item{a.} $\rm{COV}_i=100\times{\sqrt{\mathrm{var}[X_i]}}/{\mathbb{E}[X_i]}$, $i=1,\ldots,20$. 
		\item{b.} Correlation coefficients among $X_{1}$--$X_{9}$ are 0.5.
		\item{c.} $\rm{COV}_i$ for $i=19,20$ does not exist in its standard definition due to the \emph{zero} mean.		
	\end{tablenotes}			
	\label{table4}
\end{table*}
%
\subsubsection{Quantity of interest} \label{sec:5.4:2}

\begin{figure*}
	\begin{center}
		\includegraphics[angle=0,scale=0.85,clip]{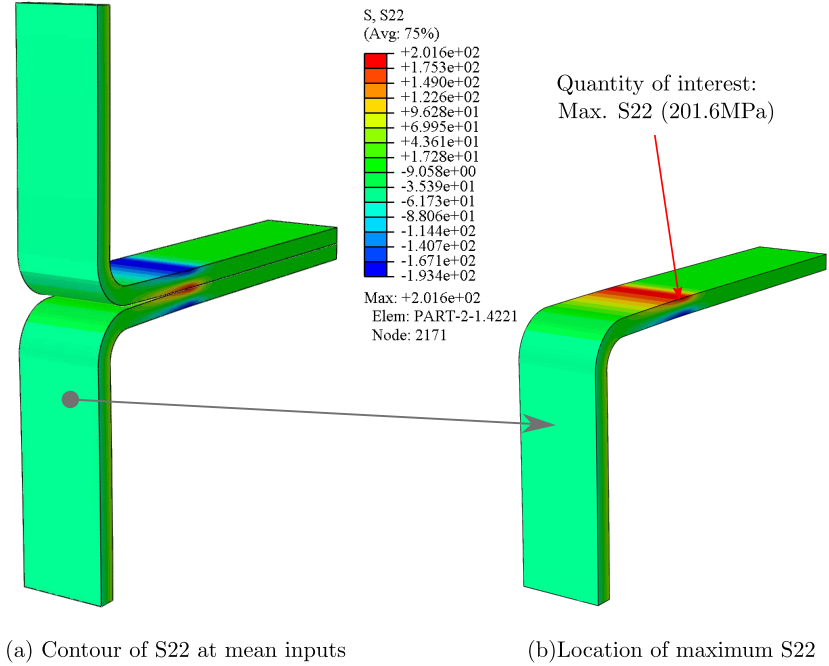}
	\end{center}
	\caption{Deterministic FEA results for the C/SiC composite T-joint. Figures~a and b show the stress in the 2-direction (S22) when the peeling displacement $U$ is applied to the T-joint and when setting the inputs to their mean values.}
	\label{fig7}
\end{figure*}

The peel test determines the peeling behavior of the C/SiC joint (i.e., the crack initiation, propagation, and arrest on the interface layer between two C/SiC substrates). Figures~\ref{fig7}a and \ref{fig7}b show the contour plots of the resultant stress (S22) in the $2$-direction when the peeling displacement load $U$ is applied to the T-joint.  Figure~\ref{fig7}b also shows the maximum value of S22. The location of the maximum value in the S22 contour can vary from sample to sample due to the randomness of the inputs $\bX$. We consider the maximum value of S22 as the output of interest for $\rm{CVaR}$ estimation since the maximum S22 can be critical for the damage initiation of the C/SiC composite materials during the peeling behavior.

\subsubsection{High-fidelity and low-fidelity model} \label{sec:5.4:3}

\begin{figure*}
	\begin{center}
		\includegraphics[angle=0,scale=0.85,clip]{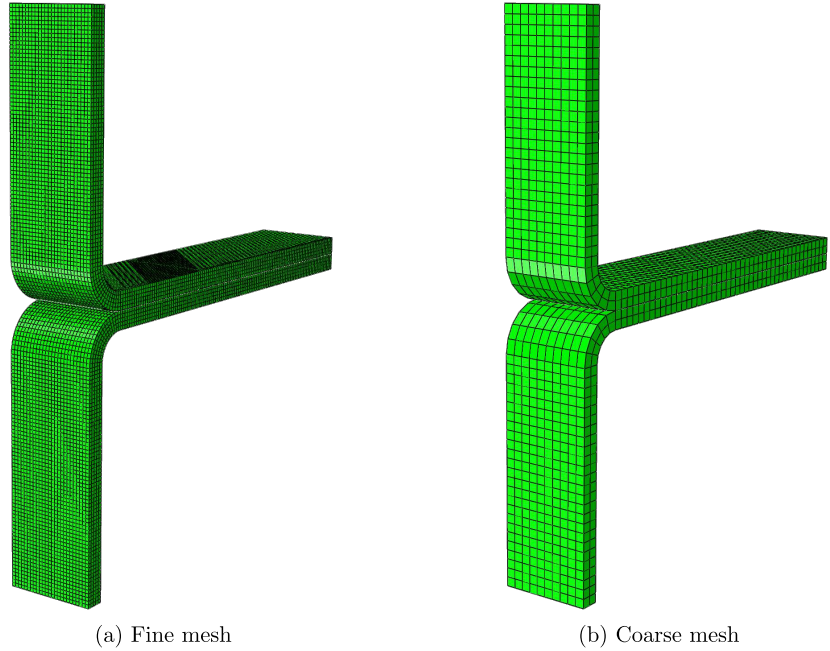}
	\end{center}
	\caption{The finite element mesh for the C/SiC T-joint. The fine mesh in (a) with 48,000 elements is used to generate high-fidelity output data, while the coarse mesh in (b) with 3,000 elements is used to generate low-fidelity output data.}
	\label{fig8}
\end{figure*}

The high-fidelity output is computed from a fine-meshed finite element model with 48,000 elements, and the low-fidelity model is based on a coarse mesh with only 3,000 elements. Those meshes are shown in Figures~\ref{fig8}a and b, respectively. The chosen element type is an eight-node brick element with reduced integration (C3D8R), built in ABAQUS/Explicit--version 6.14-2. Thus, the degrees of freedom (DOF) of the fine and coarse mesh models are 180,648 and 15,048, respectively. We determined the mesh size of the high-fidelity model by a convergence study with respect to the output of interest. We selected the low-fidelity model (which is hence not converged) such that it stays below a given computational budget.  

The correlation coefficient between the high-fidelity model $y_H$ and the low-fidelity model $y_L$ is $\rho_{H,L}=0.55$, see \eqref{pcor} for its computation. This indicates that the low-fidelity and high-fidelity model are still fairly correlated. Thus, we use the low-fidelity model to compute the DD-GPCE-Kriging surrogate, which is then used to learn the biasing density for MFIS. 

\subsubsection{Results} \label{sec:5.3:4}
\begin{table*}
	\caption{$\mathrm{CVaR}_{\beta}$ estimates ($\beta=0.99$) of the maximum value of S22 in the C/SiC T-joint via the DD-GPCE-Kriging and DD-GPCE-based MCS methods, DD-GPCE-Kriging-based MFIS, and standard MCS.}
	\vspace{-0.1in}
	\begin{spacing}{1.2}
		\begin{center}
			\small
			\resizebox{\textwidth}{!}{
				\begin{tabular}{clcccccccl}
					\toprule 
					\multicolumn{3}{c}{} & \multicolumn{3}{l}{Maximum value of S22 (MPa)} & \multicolumn{3}{l}{Model evaluations} & CPU time\tabularnewline
					\multicolumn{3}{c}{Methods} & $\mathrm{CVaR}_{\beta}$ estimates$^{(\mathrm{a})}$  & MRD$^{(\mathrm{a})}$ (\%) & {{N-RMSD$^{(\mathrm{a})}$} (\%)}  & HF$^{(\mathrm{b})}$ & LF~\RNum{1}$^{(\mathrm{c})}$ & LF~\RNum{2}$^{(\mathrm{d})}$ & (hours)$^{(\mathrm{e})}$\tabularnewline
					\midrule
					\multicolumn{5}{l}{DD-GPCE-Kriging-based MCS} &  & \tabularnewline
					& \multicolumn{2}{c}{$S=1,\,m=1$} & $242.6210$ & $2.6970$ & $2.8467$ & $400$$^{(\mathrm{f})}$ & - & - & $159.80$\tabularnewline
					& \multicolumn{2}{c}{$S=1,\,m=2$} & $245.9401$ & $2.2590$ & $2.9568$ & $400$$^{(\mathrm{f})}$ & - & - & $159.80$\tabularnewline
					& \multicolumn{2}{c}{$S=1,\,m=3$} & $249.2416$ & $2.2392$ & $3.0945$ & $400$$^{(\mathrm{f})}$ & - & - & $159.80$\tabularnewline
					\multicolumn{5}{l}{DD-GPCE-based MCS} &  & \tabularnewline
					& \multicolumn{2}{c}{$S=1,\,m=1$} & $191.1512$ & $22.5908$ & $23.8108$ & $400$$^{(\mathrm{f})}$ & - & - & $159.80$\tabularnewline
					& \multicolumn{2}{c}{$S=1,\,m=2$} & $208.8489$ & $15.4653$ & $15.7966$ & $400$$^{(\mathrm{f})}$ & - & - & $159.80$\tabularnewline
					& \multicolumn{2}{c}{$S=1,\,m=3$} & $227.6264$ & $8.4541$ & $9.0896$ & $400$$^{(\mathrm{f})}$ & - & - & $159.80$\tabularnewline
					\multicolumn{5}{l}{DD-GPCE-Kriging-based MFIS } &  & \tabularnewline
					& \multicolumn{2}{l}{HF} & $246.3827$ & $0.6644$ & $0.7156$ & $400$$^{(\mathrm{g})}$ & - & $6,000$$^{(\mathrm{h})}$ & $159.80$\tabularnewline
					& \multicolumn{2}{l}{LF} & $246.9583$ & $0.9792$ & $1.2357$ & $150$$^{(\mathrm{i})}$ & $250$$^{(\mathrm{j})}$ & $6,000$$^{(\mathrm{h})}$ & $67.69$\tabularnewline
					\multicolumn{3}{l}{Standard MCS (Benchmark)} & $247.0036$ & - &  & $6,000$ & - & - & $1597.97$\tabularnewline
					\bottomrule
				\end{tabular}
			}
			\par\end{center}
	\end{spacing}
	\vspace{-0.1in}
	\begin{tablenotes}
		\scriptsize\smallskip 
		\item{a.} The estimates are averaged over $K=30$ independent trials. 
		\item{b.} The high-fidelity output is obtained by the fine mesh model in Figure~\ref{fig8}a.   
		\item{c.} The low-fidelity output is obtained by the coarse mesh model in Figure~\ref{fig8}b.   
		\item{d.} The low-fidelity output is obtained by DD-GPCE-Kriging. 			
		\item{e.} The total CPU time is (the number of FEA)$\times$(FEA run-time) in each trial. 
		\item{f.} The high-fidelity output samples are used to calculate DD-GPCE-Kriging or DD-GPCE.
		\item{g.} Two different 50\% high-fidelity output samples are used to estimate the CI-based $\epsilon$-risk region and $\rm{CVaR}_{\beta}$, respectively.	
		\item{h.} The low-fidelity output samples are used to estimate the CI-based $\epsilon$-risk region.
		\item{i.} The high-fidelity output samples are used to estimate the $\rm{CVaR}_{\beta}$.
		\item{j.} The low-fidelity output samples are used to calculate DD-GPCE-Kriging. %
	\end{tablenotes}
	\label{table5}
\end{table*}
Table~\ref{table5} summarizes the $\rm{CVaR}$ estimates, when $\beta=0.99$, by the DD-GPCE-Kriging and DD-GPCE-based MCS methods and DD-GPCE-Kriging-based MFIS with two options (Option ``LF" trains the DD-GPCE-Kriging with the low-fidelity model and Option ``HF" with the high-fidelity model). Compared to the benchmark CVaR estimate by the standard MCS of $6,000$ high-fidelity output samples, the univariate ($S=1$) DD-GPCE-Kriging methods of the first through third-order ($m=1-3$) approximations have  a MRD of 2.2--2.7\%. While the MRD values of both DD-GPCE-Kriging and DD-GPCE decrease as the number of $m$ increases, DD-GPCE-Kriging-based estimates are more accurate than those of the DD-GPCE method. It takes $159.80$ CPU hours to obtain CVaR estimates by the DD-GPCE or the DD-GPCE-Kriging method, which is only 10\% of the $1,597.97$ CPU hours needed for the standard MCS estimate.     

The proposed MFIS-based method with HF (when the high-fidelity model is used for computing DD-GPCE-Kriging) yields the most accurate CVaR estimate, showing only 0.66\% in MRD. Moreover, the MFIS method with HF has an N-RMSD of $0.7\%$, which is the lowest among all methods in Table~\ref{table5}. It requires the same CPU time (159.80 hours) as the DD-GPCE-Kriging or DD-GPCE methods. When the proposed MFIS-based method is integrated with LF, as expected, the required CPU time decreases to 67.69 hours, which is only 42\% of the $159.80$ hours by the HF option. This is because the DD-GPCE-Kriging is calculated from cheap low-fidelity model evaluations in the LF option. The DD-GPCE-Kriging-based MFIS with the LF option yields a CVaR estimate with $0.98$\% in MRD, which is more accurate when compared with the MCS strategy that uses DD-GPCE-Kriging to sample from. This example illustrates that the proposed MFIS-based method is not only accurate but also computationally efficient when the LF option is used.   
\subsection{Guidelines for choosing a CVaR estimation strategy}  
We summarize our findings that are supported by the numerical results in this paper. The intent is to provide guidelines for the use of surrogates in risk estimation. We categorize ``efficiency" into low, medium, and high based on the computational cost that it takes to get accurate estimates. The standard MCS has low efficiency due to the prohibitive cost of evaluating the high-fidelity model. 
\paragraph{\textbf{Using the DD-GPCE-Kriging surrogate for Monte-Carlo simulation}}~\\
The method is appropriate for treating moderately to highly nonlinear and smooth output functions. However, it only has medium efficiency. The accuracy of CVaR estimates highly relies on the precision of the DD-GPCE-Kriging model. Thus, the MCS method is limited in providing accurate estimates for a nonsmooth output function.

\paragraph{\textbf{Using the DD-GPCE-Kriging surrogate for multifidelity importance sampling}}
\begin{itemize} 
	\item[(a)] Using high-fidelity evaluations to learn the surrogate (`HF'). The method is appropriate for treating highly nonlinear and nonsmooth output functions, and its efficiency is medium. The MFIS method provides unbiased CVaR estimates obtained by high-fidelity output samples.
	\item[(b)] Using low-fidelity evaluations to learn the surrogate (`LF'). The method is appropriate for treating highly nonlinear and nonsmooth output functions, and its efficiency is high. The LF option determines DD-GPCE-Kriging using low-fidelity output samples to further speed up the process of the MFIS.
\end{itemize} 
Lastly, we note that both proposed methods outperformed the MCS-based CVaR estimation when the DD-GPCE model is used for sampling. Consequently, the proposed DD-GPCE-Kriging model is a superior surrogate for uncertainty quantification, and, specifically, risk estimation.  
The DD-GPCE-Kriging-based MFIS method is especially suited to highly nonlinear and nonsmooth output responses. Moreover, the proposed MFIS-based method can achieve higher efficiency when the DD-GPCE-Kriging surrogate is computed with the low-fidelity surrogate.

\section{Conclusions} \label{sec:6} 
We presented novel computational methods for Conditional Value-at-Risk (CVaR) estimation for nonlinear systems under high-dimensional and dependent inputs. 
First, we proposed a new surrogate model, \textit{DD-GPCE-Kriging}, as a fusion of DD-GPCE and Kriging to approximate a highly nonlinear and nonsmooth random output when the inputs are high-dimensional and dependent. We integrated DD-GPCE-Kriging with two sampling-based CVaR estimation methods: standard Monte Carlo Sampling (MCS) and Multifidelity Importance Sampling (MFIS).
The proposed MCS-based method samples from the computationally efficient DD-GPCE-Kriging surrogate and is shown to be accurate in the presence of high-dimensional and dependent random inputs. However, the proposed method may produce a biased CVaR estimate as it relies on the approximation quality of the DD-GPCE-Kriging. Therefore, we leveraged the MFIS method and used the newly-constructed DD-GPCE-Kriging surrogate to determine the biasing density for importance sampling. The importance-sampling-based CVaR estimate uses the high-fidelity model only. To further speed up the MFIS process, we computed a second DD-GPCE-Kriging surrogate from lower-fidelity model simulations, which in all numerical results, still provided a good biasing density. 

The numerical results for two mathematical functions with bivariate random inputs demonstrated that the proposed DD-GPCE-Kriging method is better suited to handle a highly nonlinear or nonsmooth output than the DD-GPCE method. The example also shows that the proposed MFIS-based method provides more accurate CVaR estimates than the MCS-based method when the quantity of interest is nonsmooth.   

We evaluated the proposed CVaR estimation methods in aerospace and defense engineering applications focusing on composite structures. 
The scalability of the proposed methods and their applicability to complex engineering problems were demonstrated by solving the two-dimensional glass/vinylester laminate (lightweight and high stiffness) composite problem, including 28 (partly dependent) random inputs and the three-dimensional C/SiC (lightweight and heat resistant) composite T-joint problem, including 20 (partly dependent) random inputs. In the glass/vinylester composite laminate problem, the proposed MFIS-based method achieves a speedup of $104$x compared to standard MCS using the high-fidelity model, while producing a CVaR estimate with $1.15$\% error.
\section*{Acknowledgment} \label{sec:7} 

\noindent\textbf{\textit{Funding}} This material is based on research sponsored by the Air Force Research Lab (AFRL) and Defense Advanced Research Projects Agency (DARPA) under agreement number FA8650-21-2-7126. The U.S. Government is authorized to reproduce and distribute reprints for Governmental purposes notwithstanding any copyright notation thereon.
\newline\newline
\noindent\textbf{\textit{Disclaimer}} The views and conclusions contained herein are those of the authors and should not be interpreted as necessarily representing the official policies or endorsements, either expressed or implied, of the AFRL and DARPA or the U.S. Government.

\bibliography{ref}
\bibliographystyle{abbrv}

\begin{appendices}
	\section{Generalized polynomial chaos expansion\label{sec:appx1}}
	%
	When $\bX=(X_1,\ldots,X_N)^\intercal$ consists of statistically dependent random variables, the resultant probability measure is in general not of a product type, meaning that the joint distribution of $\bX$ cannot be obtained strictly from its marginal distributions.  Consequently, measure-consistent multivariate orthonormal polynomials in $\bx=(x_1,\ldots,x_N)^\intercal$ cannot be built from an $N$-dimensional tensor product of measure-consistent univariate orthonormal polynomials.  In this case, a three-step algorithm based on a whitening transformation of the monomial basis can be used to determine multivariate orthonormal polynomials consistent with an arbitrary, non-product-type probability measure $f_{\bX}(\bx)\rd\bx$ of $\bX$, see \cite{lee2020practical}. 
	
	Let $\bj:=(j_1,\ldots,j_N) \in \nat_0^N$ be an $N$-dimensional multi-index.  For a realization $\bx=(x_1,\ldots,x_N)^\intercal \in \mathbb{A}^N \subseteq \real^N$ of $\bX$, a monomial in the real variables $x_1,\ldots,x_N$ is the product $\bx^{\bj}=x_1^{j_1}\ldots x_N^{j_N}$ with a total degree $|\bj|=j_1+\cdots+j_N$.  Consider for each $m \in \nat_0$ the elements of the multi-index set 
	\[
	\cJ_m:=\{ \bj \in \nat_0^N: |\bj|\le m\},
	\]
	which is arranged as $\bj^{(1)},\ldots,\bj^{(L_{N,m})}$, $\bj^{(1)}=\textbf{0}$, according to a monomial order of choice.  The set $\cJ_m$ has cardinality $L_{N,m}$ obtained as
	\begin{align}
		L_{N,m}:=|\cJ_m|=\sum_{l=0}^m \binom{N+l-1}{l}=\binom{N+m}{m}.
		\label{A1}
	\end{align}
	Let us denote by
	\begin{align}
		{\mathbf{\Psi}}_{m}(\bx)=({\Psi}_1(\bx),\ldots,{\Psi}_{L_{N,m}}
		(\bx))^{\intercal}
	\end{align}
	an $L_{N,m}$-dimensional vector of multivariate orthonormal polynomials that is consistent with the probability measure $f_{\bX}(\bx)\rd{\bx}$ of $\bX$. Consequently, any output random variable $y(\bX)\in L^2(\Omega, \cF, \mathbb{P})$ can be approximated by the $m$th-order GPCE\footnote{The GPCE in \eqref{gpce} should not be confused with that of \citep*{xiu02}. The GPCE, presented here, is meant for an arbitrary dependent probability distribution of random input. In contrast, the existing PCE, whether classical \citep*{wiener38} or generalized \citep*{xiu02}, still requires independent random inputs.}
	\begin{align}
		y_m(\bX)=\displaystyle\sum_{i=1}^{L_{N,m}}c_i\Psi_i(\bX)
		\label{gpce}
	\end{align}
	of $y(\bX)$, comprising $L_{N,m}$ basis functions with expansion coefficients 
	\begin{align}
		c_i:=\displaystyle\int_{\mathbb{A}^N}y(\bx)\Psi_i(\bx)f_{\bX}(\bx)\rd\bx, i=1,\ldots,L_{N,m}.
	\end{align} 
	Here, the orthonormal polynomials $\Psi_i(\bX)$, $i=1,\ldots,L_{N,m}$, are determined numerically via three steps, which we implement by replacing $\cJ_{N,S,m}$ with  $\cJ_m$, see the following Appendix~\ref{sec:appx2}. We refer to \citep*{lee2020practical} for more details.		
	In this work, we refer to the GPCE in \eqref{gpce} as regular GPCE to distinguish it from the DD-GPCE.

	\section{Three steps for measure-consistent orthonormal polynomials for DD-GPCE\label{sec:appx2}}
	For $N\in\nat$, denote by $\{1,\ldots,N\}$ an index set and $\cU\subseteq\{1,\ldots,N\}$ a subset (including the empty set $\emptyset$) with cardinality $0\leq|\cU|\leq N$.  The complementary subset of $\cU$ is denoted by $\cU^{c}:=\{1,\ldots,N\}\backslash \cU$.  For each $m\in\nat_0$ and $0\leq S \leq N$, we define the reduced multi-index set
	\begin{equation} 
		\begin{array}{rcl}
			\cJ_{S,m}&:=&\left\{\bj=(\bj_{\cU},\textbf{0}_{\cU^c})\in\nat_0^N:\bj_{\cU}\in\nat^{|\cU|},~|\cU|\leq|\bj_{\cU}|\leq m,\right.\\&&\left.0\leq |\cU| \leq S\right\},
		\end{array} 
	\end{equation} 
	which is arranged as $\bj^{(1)},\ldots,\bj^{(L_{N,S,m})}$,  $\bj^{(1)}=\textbf{0}$, according to a monomial order of choice and where $|\bj_{\cU}| := j_{i_1}+\cdots+j_{i_{|\cU|}}$.  Here, $(\bj_{\cU},\textbf{0}_{\cU^c})$ denotes an $N$-dimensional multi-index whose $i$th component is $j_i$ if $i\in \cU$ and $0$ if $i\notin \cU$.
	The set $\cJ_{S,m}$ represents a subset of $\cJ_m$ determined from the chosen $S$, where only at most $S$-variate basis functions are preserved, that are relevant for the $S$th-variate DD-GPCE approximation. As a result, we have that
	\[
	L_{N,S,m} \leq L_{N,m},
	\]
	i.e., the DD-GPCE never has more terms than the regular GPCE; in most cases, it will have significantly less terms.
	For $\bx=(x_1,\ldots,x_N)^{\intercal}\in\mathbb{A}^N\subseteq\real^N$ we then define the basis vector for the DD-GPCE as
	\[
	\mathbf{\Psi}_{S,m}(\bx):=(\Psi_i(\bx),\ldots,\Psi_{L_{N,S,m}}(\bx))^{\intercal},
	\]
	which is an $L_{N,S,m}$-dimensional vector of multivariate orthonormal polynomials that is consistent with the probability measure $f_{\bX}(\bx)\rd\bx$ of $\bx$. The orthonormal polynomials are determined by the following three steps.  
	
	\begin{enumerate}
		[labelwidth=1.4cm,labelindent=5pt,leftmargin=1.6cm,label=\bfseries Step  \arabic*.,align=left]
		\setcounter{enumi}{0}
		\item Given $0\leq S \leq N$ and $S \leq m < \infty$, create an $L_{N,S,m}$-dimensional column vector
		\begin{equation}	 
			\bM_{S,m}(\bx)=(\bx^{\bj^{(1)}},\ldots,\bx^{\bj^{(L_{N,S,m})}})^{\intercal}
			\label{appx2:1}
		\end{equation}
		
		of monomials whose elements are the monomials $\bx^{\bj}$ for $\bj\in\cJ_{S,m}$ arranged in the aforementioned order.  For $\cU\subseteq \{1,\ldots,N\}$, let $\bx_{\cU}:=(x_{i_1},\ldots,x_{i_{|\cU|}})^{\intercal}$, $1 \leq i_1 < \cdots < i_{|\cU|} \leq N$, be a subvector of $\bx$.  The complementary subvector is defined by $\bx_{\cU^c}:=\bx_{\{1,\ldots,N\}\backslash \cU}$.  Then, for $\bj\in\cJ_{S,m}$,
		\[
		\bx^{\bj}=\bx_{\cU}^{\bj_{\cU}}{\textbf{0}_{\cU^c}}^{\bj_{\cU^c}}=\bx_{\cU}^{\bj_{\cU}}.
		\]  
		Hence, $\bM_{S,m}(\bx)$ is the monomial vector in $\bx_{\cU}=(x_{i_1},\ldots,x_{i_{|\cU|}})^{\intercal}$ of degree $0\leq |\cU| \leq S$ and $|\cU| \leq |\bj_{\cU}| \leq m$. 
		\item
		Construct an $L_{N,S,m} \times L_{N,S,m}$ monomial moment matrix of $\bM_{S,m}(\bX)$, defined as
		\begin{equation}
			\begin{split}
				{\bG}_{S,m}&:= \Exp[\bM_{S,m}(\bX)\bM_{S,m}^{\intercal}(\bX)]\\&:=
				\int_{\mathbb{A}^{N}}\bM_{S,m}(\bx)\bM_{S,m}^{\intercal}(\bx)
				f_{\bX}(\bx)\rd\bx.
			\end{split}
			\label{appx2:2}
		\end{equation}		
		For an arbitrary probability density $f_{\bX}(\bx)$, the matrix $\bG_{S,m}$ cannot be determined exactly, yet it can be accurately estimated with numerical integration and/or sampling methods \citep*{lee2020practical}.
		\item
		Select the $L_{N,S,m} \times L_{N,S,m}$ whitening matrix ${\bW}_{S,m}$ from the Cholesky decomposition of the symmetric positive-definite monomial moment matrix ${\bG}_{S,m}$ \citep*{rahman2018polynomial}, leading to
		\begin{equation}
			{\bG}_{S,m}={\bW}_{S,m}^{-1}{\bW}_{S,m}^{-\intercal}.
			\label{appx2:3}
		\end{equation}
		
		The whitening transformation is then used to generate multivariate orthonormal polynomials as follows:
		\begin{equation}
			{\mathbf{\Psi}}_{S,m}(\bx)={\bW}_{S,m} \bM_{S,m}(\bx).
		\end{equation}	
	\end{enumerate}

	\section{Leave-one-out cross validation estimate\label{sec:appx3}}
	From the known distribution of random inputs $\bX$ and an output function $y:\mathbb{A}^N\rightarrow\mathbb{R}$, consider an input-output data set $\{\bx^{(l)},y(\bx^{(l)})\}_{l=1}^{L'}$ of number $L'\in\mathbb{N}$. Given an autocorrelation function $R(|\bx-\bx'|;\btheta)$ between two input realizations $\bx$ and $\bx'$, the hyper-parameters $\btheta=(\theta_1,\ldots,\theta_N)^{\intercal}$ can be determined by the leave-one-out cross-validation estimate method \cite{bachoc2013cross} which solves   
	\begin{align}\label{ac:1} 
		\argmin_{\btheta\in\mathbb{R}^N}\big[\bb^{\intercal}\mathbf{R}^{-1}\diagmat{\mathbf{R}^{-1}}^{-2}\mathbf{R}^{-1}\bb\big],
	\end{align}
	where 
	\begin{equation}
		\begin{split}
			\bR &:=
			\begin{bmatrix}
				R(\bx^{(1)},\bx^{(1)};\btheta) & \cdots &  R(\bx^{(1)},\bx^{(L')};\btheta) \\
				\vdots                   & \ddots &  \vdots                           \\
				R(\bx^{(L')},\bx^{(1)};\btheta) & \cdots &  R(\bx^{(L')},\bx^{(L')};\btheta) 
			\end{bmatrix}\\
		\end{split}
	\end{equation}
	and $\bb=(y(\bx^{(1)}),\ldots,y(\bx^{(L')}))^{\intercal}$ is an $L'$-dimensional column vector of output evaluated at each input $\bx^{(l)}$, $l=1,\ldots,L'$.
	
	The optimal value of $\btheta$ from \eqref{ac:1} is used to determine the mean $\bar{y}_{S,m}(\bx)$ and variance $\bar{\sigma}_{S,m}^2(\bx)$ of the $S$-variate, $m$th-order DD-GPCE-Kriging, as explained in Section~\ref{sec:3.2}.     
\end{appendices}
\end{document}

%% file: ListPackages.tex
\usepackage[margin = 1.2in]{geometry}
\usepackage[T1]{fontenc}
\usepackage[english]{babel}
\usepackage[latin1]{inputenc}
\usepackage{mathtools}
\usepackage{amsmath}
\usepackage{amsfonts}
\usepackage{amsthm}
\usepackage{amssymb}
\usepackage{array}
\usepackage{subcaption}
\usepackage{siunitx}

\usepackage{color}
\usepackage{url}
\usepackage{cleveref}
\usepackage{graphicx}
\usepackage{breakcites}
\usepackage{soul} 
\usepackage{enumitem}
\usepackage[title,titletoc,toc]{appendix}



%% file: mymacros.tex
\newtheorem{remark}{Remark}

\newtheorem{problem}{Problem}

{\noindent {\textbf{Proof}:} }%
{\hfill $\Box$ \\[1ex] }

\newcommand{\bit}{\begin{itemize}}
\newcommand{\eit}{\end{itemize}}
\newcommand{\ben}{\begin{enumerate}}
\newcommand{\een}{\end{enumerate}}

\newcommand {\real} {\mathbb{R}}
\newcommand {\nat} {\mathbb{N}}

\DeclareMathOperator*{\Exp}{\mathbb{E}}

\DeclareMathOperator*{\argmin}{arg\,min}%
\newcommand{\rd}{\text{\upshape d}} 

\newcommand{\diagmat}[1]{\ensuremath{{\mathrm{diag}}\left(#1\right)}}

\newcommand{\tk}{\ensuremath{\tilde{k}}}

\newcommand{\bA}{\ensuremath{\mathbf{A}}}

\newcommand{\bG}{\ensuremath{\mathbf{G}}}

\newcommand{\bM}{\ensuremath{\mathbf{M}}}

\newcommand{\bR}{\ensuremath{\mathbf{R}}}

\newcommand{\bW}{\ensuremath{\mathbf{W}}}
\newcommand{\bX}{\ensuremath{\mathbf{X}}}

\newcommand{\bZ}{\ensuremath{\mathbf{Z}}}

\newcommand{\bb}{\ensuremath{\mathbf{b}}}
\newcommand{\bc}{\ensuremath{\mathbf{c}}}

\newcommand{\bj}{\ensuremath{\mathbf{j}}}

\newcommand{\br}{\ensuremath{\mathbf{r}}}

\newcommand{\bx}{\ensuremath{\mathbf{x}}}

\newcommand{\bz}{\ensuremath{\mathbf{z}}}

\newcommand {\bPsi} {\mbox{\boldmath $\Psi$}}

\newcommand {\btheta} {\mbox{\boldmath $\theta$}}

\newcommand{\cB}{\ensuremath{\mathcal{B}}}

\newcommand{\cD}{\ensuremath{\mathcal{D}}}

\newcommand{\cF}{\ensuremath{\mathcal{F}}}
\newcommand{\cG}{\ensuremath{\mathcal{G}}}

\newcommand{\cJ}{\ensuremath{\mathcal{J}}}

\newcommand{\cU}{\ensuremath{\mathcal{U}}}

